\documentclass{amsart}
\usepackage[english]{babel}
\usepackage[cp1250]{inputenc}
\usepackage{amsmath,amssymb}
\usepackage{microtype}
\usepackage{float}
\usepackage{enumitem}
\usepackage{mathrsfs}

\usepackage{tikz}
\tikzset{font={\fontsize{10pt}{12}\selectfont}}
\usetikzlibrary{matrix,arrows}

\newtheorem{thm}{Theorem}[section]
\newtheorem{prop}[thm]{Proposition}
\newtheorem{lem}[thm]{Lemma}
\newtheorem{cor}[thm]{Corollary}
\newtheorem{claim}[thm]{Claim}
\newtheorem*{claim*}{Claim}

\theoremstyle{definition}

\newtheorem{remark}[thm]{Remark}

\newtheorem{notation}[thm]{Notation}
\newtheorem{alg}[thm]{Algorithm}

\def\NN{\mathbb N}
\def\ZZ{\mathbb Z}

\def\QQ{\mathbb Q}
\def\CP{\mathbb CP}

\DeclareMathOperator{\spinc}{Spin^c}

\DeclareMathOperator{\tb}{tb}
\DeclareMathOperator{\rot}{rot}

\DeclareMathOperator{\PD}{PD}

\begin{document}

\title[Tight structures on small Seifert fibered $L$-spaces]{Classification of tight contact structures on small Seifert fibered $L$-spaces}

\author{Irena Matkovi\v{c}}
\address{Department of Mathematics, Central European University, 1051 Budapest, Hungary}
\email{matkovic\_irena@phd.ceu.edu}

\begin{abstract}
The Ozsv\'ath-Szab\'o contact invariant is a complete classification invariant for tight contact structures on small Seifert fibered $3$-manifolds which are $L$-spaces.
\end{abstract}


\maketitle


\section{Introduction}

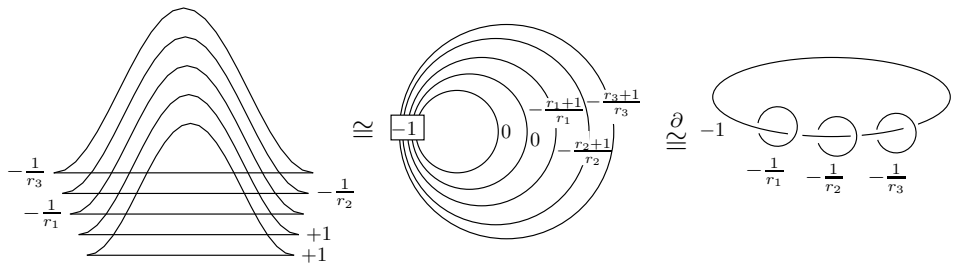
\begin{figure}[H]
\resizebox{\textwidth}{!}{%
\begin{tikzpicture}
\begin{scope}[shift={(-2,0)}]
\begin{scope}[shift={(0,0)}]
\draw [domain=-0.5*pi:1.5*pi, scale=0.5] plot (\x, {2*sin(\x r)});
\draw [scale=0.5] (-0.5*pi,-2)--(1.5*pi,-2);
\draw [scale=0.5] (-0.5*pi,-2) circle (0.1pt) node[left,scale=0.77]{$-\frac{1}{r_3}$};
\end{scope}
\begin{scope}[shift={(.06,-.3)},scale=0.95]
\draw [domain=-0.5*pi:1.5*pi, scale=0.5] plot (\x, {2*sin(\x r)});
\draw [scale=0.5] (-0.5*pi,-2)--(1.5*pi,-2) node[right,scale=0.77]{$-\frac{1}{r_2}$};
\end{scope}
\begin{scope}[shift={(.12,-.6)},scale=0.9]
\draw [domain=-0.5*pi:1.5*pi, scale=0.5] plot (\x, {2*sin(\x r)});
\draw [scale=0.5] (-0.5*pi,-2)--(1.5*pi,-2);
\draw [scale=0.5] (-0.5*pi,-2) circle (0.1pt) node[left,scale=0.77]{$-\frac{1}{r_1}$};
\end{scope}
\begin{scope}[shift={(.18,-.9)},scale=0.85]
\draw [domain=-0.5*pi:1.5*pi, scale=0.5] plot (\x, {2*sin(\x r)});
\draw [scale=0.5] (-0.5*pi,-2)--(1.5*pi,-2)node[right,scale=0.7]{$+1$};
\end{scope}
\begin{scope}[shift={(.24,-1.2)},scale=0.8]
\draw [domain=-0.5*pi:1.5*pi, scale=0.5] plot (\x, {2*sin(\x r)});
\draw [scale=0.5] (-0.5*pi,-2)--(1.5*pi,-2)node[right,scale=0.7]{$+1$};
\end{scope}
\end{scope}

\draw[white] (.7,-0.5) circle (0pt) node[right, black]{$\cong$};

\begin{scope}[shift={(-0.8,0)}]
\draw (3.5,-0.5) circle (1.3cm);
\draw [fill=white,white] (4.5,-0.4) rectangle (4.9,0);
\draw (4.8,-0.4) circle (0pt) node[above,scale=0.7]{$-\frac{r_3+1}{r_3}$};
\draw (3.37,-0.5) circle (1.13cm);
\draw [fill=white,white] (4.2,-0.9) rectangle (4.5,-0.5);
\draw (4.45,-0.5) circle (0pt) node[below,scale=0.7]{$-\frac{r_2+1}{r_2}$};
\draw (3.2,-0.5) circle (0.9cm);
\draw [fill=white,white] (3.8,-0.45) rectangle (4.2,-0.1);
\draw (4.1,-0.5) circle (0pt) node[above,scale=0.7]{$-\frac{r_1+1}{r_1}$};
\draw (3.05,-0.5) circle (0.7cm);
\draw (3.7,-0.6) circle (0pt) node[right,scale=0.7]{$0$};
\draw (2.9,-0.5) circle (0.5cm);
\draw (3.35,-0.5) circle (0pt) node[right,scale=0.7]{$0$};
\draw [fill=white] (2.1,-0.6) rectangle (2.5,-0.3) node[below left,scale=0.7]{$-1$};
\end{scope}

\draw[white] (4.5,-0.5) circle (0pt) node[right,black]{$\stackrel{\partial}{\cong}$};

\begin{scope}[scale=0.8, shift={(-.7,-.1)}]
\draw (7.5,-0.5) circle (0pt) node[left,scale=0.7]{$-1$};
\draw (9,0) ellipse (1.8cm and 0.6cm);
\draw [fill=white,white] (8.35,-0.45) rectangle (8.55,-0.65);
\draw [fill=white,white] (9.3,-0.6) rectangle (9.45,-0.4);
\draw [fill=white,white] (10.1,-0.35) rectangle (10.3,-0.55);
\draw (7.9,-.55) arc (-160:170:0.3);
\draw[white] (8,-0.8) circle (0pt) node[below,black,scale=0.77]{$-\frac{1}{r_1}$};
\draw (8.8,-.7) arc (-160:160:0.3);
\draw[white] (8.9,-.9) circle (0pt) node[below,black,scale=0.77]{$-\frac{1}{r_2}$};
\draw (9.7,-.65) arc (-160:170:0.3);
\draw[white] (9.85,-.9) circle (0pt) node[below,black,scale=0.77]{$-\frac{1}{r_3}$};
\end{scope}

\end{tikzpicture}}
\caption{Contact structures on $M(-1;r_1,r_2,r_3)$, followed by the smoothened surgery diagram of the underlying $3$-manifold and its standard presentation; when referring to them as $4$-manifolds, we assume inverse slam-dunks to be done.}
\label{fig:SFS}
\end{figure}

By small Seifert fibered 3-manifold we refer to Seifert fibration over the sphere $S^2$ with three singular fibers, standardly given as $M(e_0;r_1,r_2,r_3)$ where $e_0\in\ZZ$ and $r_i\in\QQ\cap(0,1)$ with $r_1\geq r_2\geq r_3  $. For a surgery presentation of this manifold, see the right diagram of Figure \ref{fig:SFS}.

$L$-spaces (by definition, Heegaard Floer homology lens spaces) among Seifert fibered manifolds are geometrically characterized by absence of transverse contact structures \cite[Theorem 1.1]{LS.III}, which is essential for our classification. The restriction can be simply described in terms of the Seifert constants: $L$-spaces are all manifolds with $e_0\geq 0$ and with $e_0\leq -3$, while for $e_0=-1,-2$ some explicit numerical inequalities (see Subsection \ref{Ss4}) are imposed on the triple $(r_1,r_2,r_3)$.

Problem to classify tight contact structures up to contact isotopy is usually asked for prime atoroidal manifolds; the first because tight contact structures respect connected sum decomposition of $3$-manifolds, the second because an embedded essential torus is a known source of infinitely many non-isotopic tight structures. Small Seifert fibered manifolds, beside hyperbolic ones, share these properties. On the other hand, the existence question for Seifert manifolds has been completely answered by Lisca-Stipsicz \cite{LS}: only the ones which belong to a one-parameter family of $(2n-1)$-surgeries on the torus knot $T_{2,2n+1}$ (equivalently, which are orientation preserving diffeomorphic to $M(-1;\frac{1}{2},\frac{n}{2n+1},\frac{1}{2n+3})$ for some $n\in\NN$) do not admit any tight structure. Classification then arises from the comparison of bounds: the lower bound is obtained constructively by contact surgery complemented with the use of invariants, and for the upper bound the convex surface theory is applied.

The main invariant in the classification of tight Seifert fibered manifolds is the maximal twisting number (that is, the difference between the contact framing and the fibration framing, maximized in the smooth isotopy class of a regular fiber) -- applied in convex surface theory, it allows one to give upper bounds on the number of tight structures. By the results of Wu \cite{Wu} all tight contact structures when $e_0\leq -2$ have negative maximal twisting, while for $e_0\geq 0$ they are all zero-twisting; in the work of Ghiggini \cite{G} the negative maximal twisting is further related to the existence of transverse contact structures. This, in the case of $L$-spaces,  results in a simple division: maximal twisting is equal to zero when $e_0\geq -1$, and has value $-1$ when $e_0\leq -2$. The fixed maximal twisting of a regular fiber in all the cases gives some unique contact structure on the complement of singular fibers relative to boundary, pushing the classification into tubular neighborhoods of the three singular fibers. 

The classification whenever $e_0\neq -1$ is then finished by Legendrian surgery construction -- the diagrams are simply given by Legendrianization of standard presentation of Seifert manifold; this is done in \cite{Wu,GLS0,G} for $e_0\neq -2,-1,0,\ e_0\geq 0,\ e_0=-2$, respectively. In particular, all these tight structures are Stein fillable, and classified by the first Chern class of their fillings \cite{LM}, equivalently by their $\spinc$ structure, or closer to the present context by their contact Ozsváth-Szabó invariants \cite{P}.

The remaining case of $M(-1;r_1,r_2,r_3)$ has been partly addressed already in \cite{GLS}. Here, the constructive side needs to be attacked differently because of the existence of non-fillable tight structures. Invoking that all tight structures are zero-twisting \cite{G}, Lisca-Stipsicz gave a uniform description of all possible tight structures by certain surgery diagrams.
\begin{prop}\label{prop:SFS}\cite[Proposition 6.1]{LS.III} 
Each tight contact structure with maximal twisting equal to zero on the small Seifert fibered space $M(-1;r_1,r_2,r_3)$ is given by one of the surgery presentations of Figure \ref{fig:SFS} left.\hfill$\square$
\end{prop}
This reduces the classification problem to the recognition of tightness and isotopies between the finite collection of structures, listed by the associated Thurston-Bennequin and rotation numbers.

The underlying topological question, classification of oriented $2$-plane fields $\xi\in\Xi$ up to homotopy is given by their induced $\spinc$ structure $\mathbf t_\xi$ together with the 3-dimensional invariant $d_3(\xi)$. (Recall that $\pi_0(\Xi)$ can be identified with homotopy classes of maps $[M,S^2]$, which can be through Pontryagin-Thom construction given by framed links in $M$; here link up to oriented cobordism represents the class in $H_1(M;\ZZ)$, equivalently $\mathbf t_\xi$, while the framing corresponds to the Hopf invariant as a $3$-dimensional obstruction for homotopies between plane fields.) 

To detect tightness, we basically use the Ozsváth-Szabó contact invariant \cite{OSz.c}, implicitly expecting all tight structures to have non-vanishing one. But we address it indirectly, by showing the sufficient condition of Lisca-Stipsicz. 
\begin{thm}\label{thm:d=d3}\cite[Theorem 1.2]{LS.III} 
If for a contact structure $\xi$ of Figure \ref{fig:SFS} on Seifert fibered $L$-space $M=M(-1;r_1,r_2,r_3)$ the equality $d_3(\xi)=d(M,\mathbf t_\xi)$ holds, then its contact invariant $c(M,\xi)\in\widehat{HF}(-M,\mathbf t_\xi)$ does not vanish. \hfill$\square$
\end{thm}
Then, to confirm overtwistedness of all non-detected structures and to obtain the isotopies between tight ones, as always, convex surface theory is applied.

The observations accumulate in the confirmation of \cite[Conjecture 4.7]{S}.

\begin{thm}\label{thm}
Let $M$ be a small Seifert fibered $L$-space of the form $M(-1;r_1,r_2,r_3)$. Then a contact structure $\xi$ on $M$ is tight if and only if it is given by a contact surgery presentation of Figure \ref{fig:SFS} and its $3$-dimensional invariant $d_3(\xi)$ is equal to the $d$-invariant $d(M,\mathbf t_\xi)$. Moreover, two tight structures $\xi_1$ and $\xi_2$ on $M$ are contact isotopic if and only if their induced $\spinc$ structures $\mathbf{t}_{\xi_1}, \mathbf{t}_{\xi_2}$ are isomorphic.
\end{thm}

Expressed in terms of the Ozsv\'ath-Szab\'o contact invariant all tight structures on small Seifert fibered $L$-spaces  satisfy the following.

\begin{cor}\label{cor} 
Let $\xi$ be a contact structure on small Seifert fibered $L$-space $M=M(e_0;r_1,r_2,r_3)$. Then $\xi$ is tight if and only if its contact invariant $c(\xi)\in \widehat{HF}(-M,\mathbf t_{\xi})$ is nonzero. Moreover, two tight structures $\xi_1$ and $\xi_2$ are isotopic if and only if their contact invariants $c(\xi_1), c(\xi_2)$ coincide,  if and only if their induced $\spinc$ structures $\mathbf{t}_{\xi_1}, \mathbf{t}_{\xi_2}$ are isomorphic.
\end{cor}

\proof
If there are less than three singular fibers, the manifold considered is a lens space. Here as well as when $e_0\neq -1$, all tight structures are Stein fillable according to previous results \cite{E,Ho.I,Wu,GLS0,G}. By that and Theorem \ref{thm} above, tight structure on any considered $M$ has non-vanishing contact invariant.

The fillable structures are classified by the contact invariant due to Plamenevskaya \cite{P}. In fact, for $L$-spaces the non-trivial contact invariant of $\xi$ is the unique generator of $\widehat{HF}(-M,\mathbf t_{\xi})$, hence $\xi$ is the only tight representative of its induced $\spinc$ structure. Its $3$-dimensional invariant $d_3(\xi)$ is specified as the absolute grading of the contact invariant, which equals $d(M,\mathbf t_\xi)$. 
\endproof

Our result reduces the classification problem to a well-understood computation of invariants. Although our method does not result in the number of tight structures on a given small Seifert manifold, the problem is translated to a completely combinatorial (so not geometric) count. Indeed, in any special case the number can be easily determined by, say, a computer calculation (as here both $d_3$ and $d$ are computable, and the $\spinc$ structure can be given as an element of $H_1$). What is more, since there is a surgery presentation of considered contact manifolds, we have a very explicit description of tight structures.

\begin{remark}
In contrast to the cases with $e_0\neq -1$, not all tight structures on $M(-1;r_1,r_2,r_3)$ are fillable. Whenever $r_1+r_2<1$, in the language of \cite{LL} for manifolds of special type, the existence of Stein filling is even topologically obstructed. And as all contact structures of the form given by Figure \ref{fig:SFS} are known to be supported by open books with planar pages \cite{LS.III}, the theorem of Wendl \cite{W} implies they are not fillable at all. Most manifolds with $r_1+r_2\geq 1$ admit Stein fillable as well as non-fillable tight structures, as specified in \cite{M}.
\end{remark}

\subsection*{Overview} In Section \ref{Sec2} we explain the structure of our proof, and review main concepts behind it. Then in Section \ref{Sec3} we illustrate the suggested approach by reproving the classification on $M_p:=M(-1;\frac{1}{2},\frac{1}{2},\frac{1}{p})$ \cite{GLS}. Technical details are given in the last two sections. In Section \ref{Sec4} we establish paths of characteristic covectors separating the presentations into classes with the same contact invariant. Finally in Section \ref{Sec5}, with the help of convex surface theory, the presentations of the same class are realized to be contact isotopic, more, the failure of the tightness criterion is related to overtwistedness.

\subsection*{Acknowledgement} I am indebted to András Stipsicz for his insightful mentoring.

\section{Outline of the proof}\label{Sec2}

Following the classification scheme given in the Introduction, we need a construction, a method to detect tightness (Subsection \ref{Step1}), and finally a proof that it is complete (that is, a way to recognize overtwistedness and isotopies between possibly different presentations of the same contact structure; Subsection \ref{Step2}).

By Proposition \ref{prop:SFS}, to construct tight structures on $M(-1;r_1,r_2,r_3)$ which is an $L$-space, the contact surgery presentations of Figure \ref{fig:SFS} suffice. This gives a finite collection of contact structures, on which we need to run the following two-step analysis.

\subsection{Detect tightness}\label{Step1}
In order to detect tightness we examine the equality between the 3-dimensional homotopy invariant $d_3$ of the contact structure and the $d$-invariant of the induced $\spinc$ structure (according to Theorem \ref{thm:d=d3}). The advantage of this condition over the Ozsv\'ath-Szab\'o contact invariant is that these two invariants are easily computable. They can be described in terms of characteristic covectors on plumbings bounded by $\pm M$, which brings them into the same picture (as presented below).

We exploit two 4-manifolds, naturally arising as fillings of our $M$, one given by the smoothened surgery diagram (Figure \ref{fig:SFS} middle) -- call it $X$ -- and another given by standard smooth presentation of Seifert fibration (Figure \ref{fig:SFS} right) -- say $W=W_\Gamma$ as we can think of it as a simple plumbing along the graph $\Gamma$. They are related by $X\#\CP^2\cong W\#2\CP^2$. Further, the plumbing description of $-M$ will play the central role in what follows. Call it $W_{\Gamma'}$ where $\Gamma'$ stands for the plumbing graph, dual to $\Gamma$. 

Now, the 3-dimensional invariant of the contact bundle $\xi$ can be directly read off from the surgery presentation $d_3(\xi) =\frac{1}{4}(c^2(X,J)-3\sigma(X)-2b_2(X))+\#(+1\text{-surgeries})$ \cite[Corollary 3.6]{DGS}, where $c(X,J)$ stands for the characteristic element, determined by the $\xi$-induced almost complex structure $J$ of $X\backslash\{\text{a point in each +1-handle}\}$. While the $d$-invariant corresponds to the reversely oriented $-M$ (which bounds a negative-definite plumbing) together with an induced $\spinc$ structure $\mathbf t_\xi$. It is realized by the characteristic 2-cohomology class, which gives $\spinc$ cobordism from $S^3$ to $(-M,\mathbf t_\xi)$ whose associated map in Heegaard Floer homology decreases absolute grading the least. These are recognized by full paths \cite[Subsection 3.1]{OSz}. To establish some terminology let us recollect.

\subsubsection*{Full path}
Assume that $\Gamma'$ is a negative definite plumbing with at most one bad vertex. A $t$-tuple $(K_1,\ldots,K_t)$ of characteristic covectors on $W_{\Gamma'}$ forms a full path if its elements are connected by the following $2\PD$ steps: for some vertex $v$ with $\langle K_i,v\rangle=-v\cdot v$, the vector $K_{i+1}$ is given by $K_{i+1}=K_i+2\PD(v)$. (In particular, all its elements induce the same $\spinc$ structure on the boundary $\partial W_{\Gamma'}=-M$ as they differ only by twice generators of $H^2(W_{\Gamma'},-M;\ZZ)$. Further, their (common) degree can be computed by the formula $\frac{1}{4}(K_i^2+|\Gamma'|)$.) The path ends either by some characteristic vector $K$ which exceeds the bounds $v\cdot v\leq\langle K_i,v\rangle\leq -v\cdot v$ at some $v\in\Gamma'$ -- we will say that it drops out. Otherwise the path reaches the proper ends in the initial vector $K_1$ satisfying $v\cdot v+2\leq\langle K_1,v\rangle\leq -v\cdot v\ \text{ for all } v\in\Gamma'$ and the terminal vector $K_t$ satisfying $v\cdot v\leq\langle K_t,v\rangle\leq -v\cdot v-2\ \text{ for all } v\in\Gamma'$ -- such ending full path according to \cite{OSz} determines a non-trivial element of $\widehat{HF}(-M)$.

\subsubsection*{Embedding into blown-up $\CP^2$}\hspace*{\fill} \\

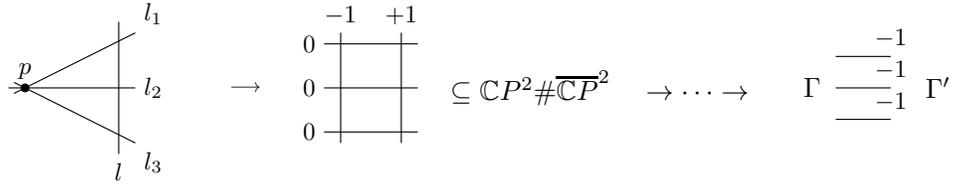
\begin{figure}[H]
\resizebox{\textwidth}{!}{%
\begin{tikzpicture}
\begin{scope}[shift={(-3.3,0)}, scale=0.7]
\filldraw[black] (0,0) circle (2pt)node[above,scale=0.9]{$p$};
\draw (-.3,0) -- (2,0)node[right,scale=0.9]{$l_2$};
\draw (-.2,.1) -- (2,-1)node[below right, scale=0.9]{$l_3$};
\draw (-.2,-.1) -- (2,1) node[above right, scale=0.9]{$l_1$};
\draw (1.7,1.2) -- (1.7,-1.2) node[below,scale=0.9]{$l$};
\end{scope}

\draw[white] (-.5,0) circle (0.1pt) node[black,scale=0.7]{$\longrightarrow$};

\begin{scope}[shift={(.5,.7)}, scale=0.7]
\draw (1.7,-.2) -- (0,-.2) node[left,scale=0.9]{$0$};
\draw (0.3,-2) -- (0.3,0);
\filldraw[black] (0.3,0) circle (0pt)node[above,scale=0.9]{$-1$};
\draw (1.4,-2) -- (1.4,0);
\filldraw[black] (1.4,0) circle (0pt)node[above,scale=0.9]{$+1$};
\draw (1.7,-1) -- (0,-1)node[left,scale=0.9]{$0$};
\draw (1.7,-1.8) -- (0,-1.8)node[left,scale=0.9]{$0$};
\end{scope}	

\draw[white] (4,0) circle (0.1pt) node[black]{$\subseteq\CP^2\#\overline{\CP}^2\ \ \ \rightarrow\cdots\rightarrow$};

\begin{scope}[shift={(7,0)}]
\draw[white] (-.3,0) circle (0.1pt) node[black]{$\Gamma$};
\draw (0,0) -- (.7,0)node[above,scale=0.9]{$-1$};
\draw (0,.4) -- (.7,.4)node[above,scale=0.9]{$-1$};
\draw (0,-.4) -- (.7,-.4)node[above,scale=0.9]{$-1$};
\draw[white] (1.3,0) circle (0.1pt) node[black]{$\Gamma'$};
\end{scope}	
\end{tikzpicture}}

\caption{Construction of the $4$-manifold $R=\CP^2\#n\overline{\CP^2}$; at the end, the two plumbings are glued together along the exceptional spheres from the last blow-up of each singular fiber, and (not shown) all regular fibers.}
\label{fig:mntR}
\end{figure}

According to \cite[Lemma 4.2]{LS} we can embed $M$ as a hypersurface in a closed oriented 4-manifold $R$ so that $R\backslash \nu(M)=W_\Gamma\cup W_{\Gamma'}$. The configuration of both intersection graphs $\Gamma,\Gamma'$ is obtained by blowing-up the initial lines $l_1,l_2,l_3\subset\CP^2: l_1\cap l_2\cap l_3=\{p\}$, and $l\subset\CP^2: p\notin l$ (see Figure \ref{fig:mntR}).

Denote standard generators of $H_2(R;\ZZ)$ as follows: $h\ (h^2=1)$ for the initial $\CP^1\subset\CP^2$, and $e_i\ (e_i^2=-1)$ for exceptional curves. Then the above description of the embedding $W_\Gamma\cup W_{\Gamma'}\hookrightarrow R$ gives:
\begin{itemize}[leftmargin=.7cm]
\item $\{z=\mbox{ center of }\Gamma\}\mapsto e_1$
\item $\{z'=\mbox{ center of }\Gamma'\}\mapsto h-e_2-e_3-e_4$
\item $\{x_i=\mbox{ first vertex of the leg }L_i\subset\Gamma\}\mapsto h-e_1-e_{i+1}-\sum e_j \mbox{ for } i=1,2,3$
\item $\{v \mbox{ vertex, }v\neq z,z',x_i\}\mapsto e_j-\sum e_k$, for example\\
$\{x_i'=\mbox{ first vertex of the dual leg }L_i'\subset\Gamma'\}\mapsto e_{i+1}-\sum e_k \mbox{ for } i=1,2,3$.
\end{itemize}

We will refer to $\Gamma$ as the manifold side and $\Gamma'$ as the dual side. Throughout, we will follow the convention that primed notation belongs to the dual graph: apart from the special vertices denoted above, let $v_j^i$ be the $j^{\text{th}}$ vertex of $L_i$ and $v_j^{i'}$ the $j^{\text{th}}$ vertex of $L_i'$.

\subsubsection*{Tightness criterion}
We have $d_3$ given by a characteristic covector on $W_\Gamma$, $d$ given by a characteristic covector on $W_{\Gamma'}$, and what is more, we can glue these two plumbings along a rational homology sphere $M$, giving blown-up $\CP^2$ (called $R$). Having a characteristic covector $c$ on $R$, which agrees with the ones providing $d_3$ and $d$ on $W_\Gamma$ and $W_{\Gamma'}$ respectively, the equality $d_3(\xi)=d(M;\mathbf t_\xi)$ can be rewritten as $c^2=\sigma(R)$. This can be understood as geometrization of Theorem \ref{thm:d=d3}.
\begin{thm} \label{thm:c}\cite[Theorem 3.3]{LS}
Let $M=M(-1;r_1,r_2,r_3)$ be an $L$-space. Then a sufficient condition for a contact structure $\xi$ on $M$ to be tight (even for $c(M,\xi)\neq0$) is in the existence of a characteristic cohomology class $c\in H^2(R;\ZZ)$ such that:
\begin{enumerate}[leftmargin=.7cm]
\item[(i)] $d_3(\xi)=\frac{1}{4}((c|_{W_\Gamma})^2-3\sigma(W_\Gamma)-2b_2(W_\Gamma))+1$ ($c|_{W_\Gamma}$ corresponds to $c(X,J)$);
\item[(ii)] $c|_{W_{\Gamma'}}$ belongs to an ending full path on $W_{\Gamma'}$ (gives rise to a class of $\widehat{HF}(-M,\mathbf{t}_\xi)$);
\item[(iii)] $c^2=\sigma(R)$.\hfill$\square$
\end{enumerate}\
\end{thm}

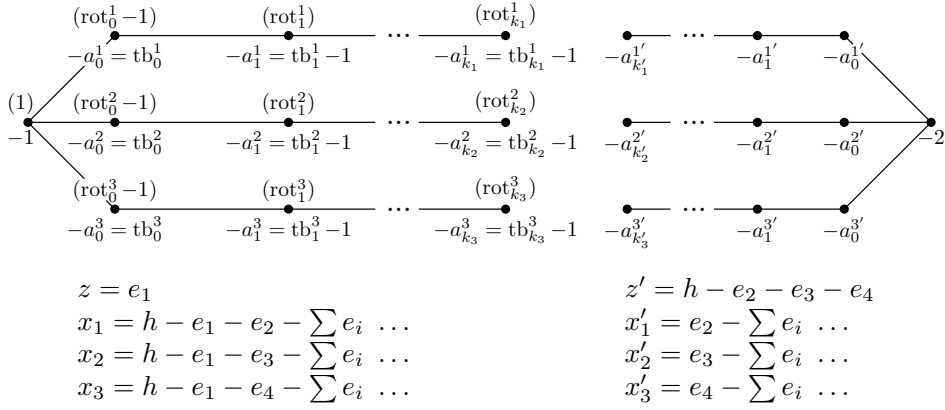
\begin{figure}[H]
\resizebox{\textwidth}{!}{%
\begin{tikzpicture}
\begin{scope}[shift={(-3.3,0)}]
\filldraw[black] (0,0) circle (1.3pt)node[above,scale=0.7]{$(1)\ \ $} node[below,scale=0.7]{$-1\ \ $};
\filldraw[black] (1,0) circle (1.3pt)node[above,scale=0.7]{$(\rot_0^2-1)$} node[below,scale=0.7]{$-a_0^2=\tb_0^2$};
\filldraw[black] (1,1) circle (1.3pt)node[above,scale=0.7]{$(\rot_0^1-1)$} node[below,scale=0.7]{$-a_0^1=\tb_0^1$};
\filldraw[black] (1,-1) circle (1.3pt)node[above,scale=0.7]{$(\rot_0^3-1)$} node[below,scale=0.7]{$-a_0^3=\tb_0^3$};
\filldraw[black] (3,1) circle (1.3pt)node[below,scale=0.7]{$-a_1^1=\tb_1^1-1$}node[above,scale=0.7]{$(\rot_1^1)$};
\filldraw[black] (3,0) circle (1.3pt)node[below,scale=0.7]{$-a_1^2=\tb_1^2-1$}node[above,scale=0.7]{$(\rot_1^2)$};
\filldraw[black] (3,-1) circle (1.3pt)node[below,scale=0.7]{$-a_1^3=\tb_1^3-1$}node[above,scale=0.7]{$(\rot_1^3)$};
\filldraw[black] (5.5,1) circle (1.3pt)node[below,scale=0.7]{$-a_{k_1}^1=\tb_{k_1}^1-1$}node[above,scale=0.7]{$(\rot_{k_1}^1)$};
\filldraw[black] (5.5,0) circle (1.3pt)node[below,scale=0.7]{$-a_{k_2}^2=\tb_{k_2}^2-1$}node[above,scale=0.7]{$(\rot_{k_2}^2)$};
\filldraw[black] (5.5,-1) circle (1.3pt)node[below,scale=0.7]{$-a_{k_3}^3=\tb_{k_3}^3-1$}node[above,scale=0.7]{$(\rot_{k_3}^3)$};
\draw (0,0) -- (1,0);
\draw (0,0) -- (.65,-.65) (.9,-.9) -- (1,-1);
\draw (0,0) -- (.65,.65) (.9,.9) -- (1,1);
\draw (1,-1) -- (4,-1)node[right]{...};
\draw (4.5,-1) -- (5.5,-1);
\draw (1,0) -- (4,0)node[right]{...};
\draw (4.5,0) -- (5.5,0);
\draw (1,1) -- (4,1)node[right]{...};
\draw (4.5,1) -- (5.5,1);
\end{scope}

\begin{scope}[shift={(3.1,0)}]
\filldraw[black] (0.5,-1) circle (1.3pt)node[below,scale=0.7]{$-a_{k_3'}^{3'}$};
\filldraw[black] (2,-1) circle (1.3pt)node[below,scale=0.7]{$-a_1^{3'}$};
\filldraw[black] (0.5,1) circle (1.3pt)node[below,scale=0.7]{$-a_{k_1'}^{1'}$};
\filldraw[black] (2,1) circle (1.3pt)node[below,scale=0.7]{$-a_1^{1'}$};
\filldraw[black] (0.5,0) circle (1.3pt)node[below,scale=0.7]{$-a_{k_2'}^{2'}$};
\filldraw[black] (2,0) circle (1.3pt)node[below,scale=0.7]{$-a_1^{2'}$};
\filldraw[black] (3,-1) circle (1.3pt)node[below,scale=0.7]{$-a_0^{3'}$};
\filldraw[black] (3,0) circle (1.3pt)node[below,scale=0.7]{$-a_0^{2'}$};
\filldraw[black] (3,1) circle (1.3pt)node[below,scale=0.7]{$-a_0^{1'}$};
\filldraw[black] (4,0) circle (1.3pt)node[below,scale=0.7]{$-2$};
\draw (3,-1) -- (4,0);
\draw (3,0) -- (4,0);
\draw (3,1) -- (3.1,.9)(3.2,.8) -- (4,0);
\draw (1.5,-1) -- (3,-1);
\draw (0.5,-1) -- (1,-1)node[right]{...};
\draw (1.5,1) -- (3,1);
\draw (0.5,1) -- (1,1)node[right]{...};
\draw (1.5,0) -- (3,0);
\draw (0.5,0) -- (1,0)node[right]{...};
\end{scope}

\begin{scope}[shift={(-3,-2.5)}]
\filldraw[white] (0,0) circle (.1pt) node[right,black,scale=.9]{$%
\begin{array}{l}z=e_1\\ x_1=h-e_1-e_2-\sum e_i\ \ldots\\ x_2=h-e_1-e_3-\sum e_i\ \ldots\\ x_3=h-e_1-e_4-\sum e_i\ \ldots\end{array}$};
\end{scope}

\begin{scope}[shift={(3.3,-2.5)}]
\filldraw[white] (0,0) circle (.1pt) node[right,black,scale=.9]{$%
\begin{array}{l}z'=h-e_2-e_3-e_4\\ x_1'=e_2-\sum e_i\ \ldots\\ x_2'=e_3-\sum e_i\ \ldots\\ x_3'=e_4-\sum e_i\ \ldots\end{array}$};
\end{scope}
\end{tikzpicture}}

\caption{Plumbing graph $\Gamma$ (left) and its dual $\Gamma'$ (right) with denoted self-intersections and evaluations of characteristic covector $c$, $(\cdot)=\langle c,v\rangle$, on the manifold side; the central and the first vertices on legs are given in generating classes of $H_2(R;\ZZ)$.}
\label{fig:graphs}
\end{figure}

For a given surgery presentation, the conditions read as follows (see Figure \ref{fig:graphs}):
\begin{enumerate}[leftmargin=.7cm]
\item[(i)] $c$ on $W_\Gamma$ is determined by $c(X,J)$, on generators of $H_2(X)$ evaluated as rotation numbers (central blow-up decreases all the (neighboring) values by $1$);
\item[(ii)] is to be checked given the constraints from (i);
\item[(iii)] we can give $c$ as $\PD(c)=\alpha h+\sum\alpha_i e_i$ where $\alpha,\alpha_i\in\{\pm1\}$, by construction fulfilling the equality $c^2=\sigma(R)$.
\end{enumerate}
 
Technically speaking, given a tuple $(c_v)_{v\in\Gamma}$ of $c$-evaluations $(c_v)=\langle c,v\rangle$ as in (i), we list all possible $\pm 1$ distributions in the expression $\PD(c)$ of (iii), and calculate corresponding values on $\Gamma'$. We know that these different $\Gamma'$-evaluations are for each given surgery presentation only different representatives of the same full path -- they all describe the same $\spinc$ structure on the boundary. We continue this path towards its ends so that we connect to it all the characteristic covectors which can be obtained by allowed $2\PD$ steps. Taking all $(c_v)$-tuples, this results in the separation of the characteristic lattice into components; denote them $\mathcal{P}_\xi$ according to the contact structure $\xi$ they belong.

In this language, Theorem \ref{thm} takes the following working form.

\begin{thm}\label{thm'}
The contact structure $\xi$ on $M(-1;r_1,r_2,r_3)$ given by surgery diagram is tight if and only if its full path $\mathcal{P}_\xi$ properly ends in the initial and terminal vector. Two such contact structures $\xi_1,\xi_2$ are isotopic if and only if their paths $\mathcal{P}_{\xi_1}, \mathcal{P}_{\xi_2}$ meet (hence, coincide).
\end{thm}

\subsection{Prove overtwistedness and describe contact isotopies}\label{Step2}
Finally, to close our classification we need that the zero elements (drop-outs) correspond to overtwistedness, and for the second part of Theorem \ref{thm} (Theorem \ref{thm'}) that elements giving the same $\widehat{HF}(M)$-generator (sharing the same path) are actually contact isotopic.

Here, convex surface theory comes in. We need to translate contact surgeries back into convex decomposition. Natural convex decomposition of the manifold $M$ separates the three singular tori from the rest of the manifold. Then the coefficients in the continued fraction expansions of the three surgeries, along with the chosen stabilizations determine basic slice decompositions of the three tori.

\subsubsection*{Contact surgery}
Contact surgery in addition to ordinary surgery prescribes for the contact structure to be preserved in the complement of a tubular neighborhood of the core link, while the extension to glued-up tori needs to be tight. The extended contact structure is  determined by the boundary slope \cite[Theorem 2.3]{Ho.I}, given by the surgery coefficient. Contact surgery diagrams \cite{DGS} encode a basic slice decomposition of the glued-up solid torus. The slope uniquely determines continued fraction blocks. Concretely, writing out the continued fraction expansion of the surgery coefficient $$-\frac{1}{r_i}=-a_0^i-\frac{1}{\ddots-\frac{1}{-a_{k_i}^i}}=[a_0^i,\ldots,a_{k_i}^i]$$ (we use the convention of \cite{GLS} with $a_j^i\geq 2$), continued fraction blocks are toric annuli in the layering of the solid torus with boundary slope $[a_{k_i}^i,\ldots,a_0^i]$, cut out by pairs of tori of slopes $[a^i_{k_i} , \ldots , a^i_{j+1} - 1]$ and $[a^i_{k_i} , \ldots , a^i_j - 1]$ (the outermost being $[a_{k_i}^i,\ldots,a_0^i]$, and the innermost $-1$). In a surgery diagram, they are represented in a chain of pushed-off knots with appropriate integral surgery coefficients (after turning into $\pm 1$ surgeries, captured by Thurston-Bennequin invariants). The remaining ambiguity in the signs of the basic slices within each continued fraction block is reflected in the choice of stabilizations of the corresponding Legendrian knot (equivalently their differences as their total number is determined by the surgery coefficient), so given by its rotation number. In the translation positive and negative basic slices in the decomposition of a continued fraction block correspond to positive and negative stabilizations (down- and up-cusps) of corresponding Legendrian knot. The loss of basic slice ordering in the transition is explained by the shuffling property \cite[Lemma 4.14]{Ho.I} of basic slices within a single block.

\medskip
What we need is to relate steps in the full path with appropriate state traversals, and drop-outs to non-tight basic slice configurations.

In other words, we have set up two ways to describe rotation numbers. Provided the surgery coefficient is fixed, they can be equally given by either the number of negative basic slices (up-cusps) or the number of negative signs on the generators forming the corresponding part of the dual leg. Then the nice thing -- to be shown -- is that the full path connections reflect the known behavior of basic slices.

\section{First example}\label{Sec3}

We illustrate our strategy on small Seifert fibered $L$-spaces $M_p:=M(-1;\frac{1}{2},\frac{1}{2},\frac{1}{p})$. The classification on these manifolds was first obtained by Ghiggini-Lisca-Stipsicz in \cite{GLS}; wherever applicable, we use their notation. First we describe tight structures on $M_p$ using Theorem \ref{thm}, then we prove Theorem \ref{thm} in this special case.

\begin{claim}\label{example}
Manifold $M_p$ admits exactly three tight contact structures up to isotopy.
\end{claim}

The finite collection of contact structures, given by Figure \ref{fig:SFS}, can be encoded in the following table of invariants:
$$
\begin{array}{ccl}
\mbox{surgery coefficient} & \ \tb\ \  & \rot (|\rot|\leq -\tb-1)\\
+1&-1&0 \\
+1&-1&0 \\
-1&-2&\rot_1\in\{-1,1\} \\
-1&-2&\rot_2\in\{-1,1\} \\
-1&-p&\rot_3\in\{-p+1, -p+3,...,p-1\}.
\end{array}
$$
As an application of the Theorem, the tightness and isotopies can be recognized solely from the induced $\spinc$ structures and the two invariants. In our case these are as follows.
$$\begin{array}{cl}
d_3(\xi)& =\frac{1}{4}(c^2(X,J)-3\sigma(X)-2b_2(X))+q\\
& =\frac{1}{4}((0,0,\rot_1,\rot_2,\rot_3)Q_X^{-1}(0,0,\rot_1,\rot_2,\rot_3)^T-3\cdot(-1)-2\cdot 5)+2.
\end{array}
$$
So, for mixed $(\rot_1,\rot_2)=(\pm1,\mp1)$, the $d_3$ is always zero, as for $(\rot_1,\rot_2)=(\pm1,\pm1)$ it runs through the values $\{\frac{2-p}{4},...,\frac{-2+3p}{4}\}$ by the step $\pm1$ as $\rot_3$ increases.

There are exactly four $\spinc$ structures  for each $p$ (as $|H_1(M_p;\ZZ)|=4$):
$$
\begin{array}{cl}H_1(-M_p;\ZZ)&=\left\langle \mu,\mu_a,\mu_b,\mu_c; 
\left(\begin{array}{cccc}1&1&1&1\\1&p&0&0\\1&0&2&0\\1&0&0&2 \end{array}\right)
\left(\begin{array}{c}\mu\\ \mu_a\\ \mu_b\\ \mu_c \end{array}\right)=0
\right\rangle\\
&=\left\{\begin{array}{l}\langle \mu_b;4\mu_b=0\rangle\cong\ZZ_4 \mbox{ for } p \mbox{ odd}\\
\langle \mu_b,\mu_c;2\mu_b=2\mu_c=0\rangle\cong\ZZ_2\oplus\ZZ_2 \mbox{ for } p \mbox{ even.} \end{array}\right.
\end{array}
$$
They can be given by the set $\{\mathbf t_1=\mathbf t_4+\mu_b, \mathbf t_2=\mathbf t_4+\mu_c, \mathbf t_3=\mathbf t_4+\mu_a, \mathbf t_4\}$. And corresponding four characteristic $2$-cohomology classes, realizing $d(-M_p,\mathbf t_i)$, are on the generators of $H_2(W_{\Gamma'})$ given by:\\
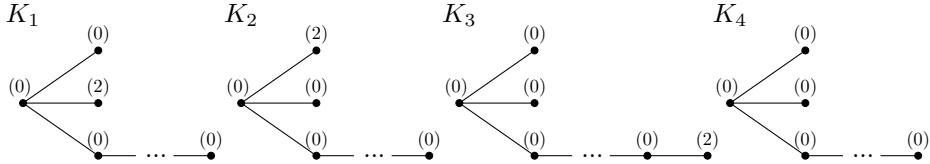
\begin{tikzpicture}
\begin{scope}[shift={(-2.9,0)}]
\filldraw[black] (0,0) circle (1.3pt)node[above,scale=0.7]{$(0)\ $};
\filldraw[black] (1,0) circle (1.3pt)node[above,scale=0.7]{$(2)$};
\filldraw[black] (1,.7) circle (1.3pt)node[above,scale=0.7]{$(0)$};
\filldraw[black] (1,-.7) circle (1.3pt)node[above,scale=0.7]{$(0)$};
\filldraw[black] (2.5,-.7) circle (1.3pt)node[above,scale=0.7]{$(0)$};
\draw (0,0) -- (1,0);
\draw (0,0) -- (1,-.7);
\draw (1,-.7) -- (1.5,-.7)node[right]{...};\draw (2,-.7) -- (2.5,-.7);
\draw (0,0) -- (1,.7);
\filldraw[white] (0,.9) circle (.1pt) node[above,black]{$K_1$};
\end{scope}

\begin{scope}[shift={(0,0)}]
\filldraw[black] (0,0) circle (1.3pt)node[above,scale=0.7]{$(0)\ $};
\filldraw[black] (1,0) circle (1.3pt)node[above,scale=0.7]{$(0)$};
\filldraw[black] (1,.7) circle (1.3pt)node[above,scale=0.7]{$(2)$};
\filldraw[black] (1,-.7) circle (1.3pt)node[above,scale=0.7]{$(0)$};
\filldraw[black] (2.5,-.7) circle (1.3pt)node[above,scale=0.7]{$(0)$};
\draw (0,0) -- (1,0);
\draw (0,0) -- (1,-.7);
\draw (1,-.7) -- (1.5,-.7)node[right]{...};\draw (2,-.7) -- (2.5,-.7);
\draw (0,0) -- (1,.7);
\filldraw[white] (0,.9) circle (.1pt) node[above,black]{$K_2$};
\end{scope}

\begin{scope}[shift={(2.9,0)}]
\filldraw[black] (0,0) circle (1.3pt)node[above,scale=0.7]{$(0)\ $};
\filldraw[black] (1,0) circle (1.3pt)node[above,scale=0.7]{$(0)$};
\filldraw[black] (1,.7) circle (1.3pt)node[above,scale=0.7]{$(0)$};
\filldraw[black] (1,-.7) circle (1.3pt)node[above,scale=0.7]{$(0)$};
\filldraw[black] (2.5,-.7) circle (1.3pt)node[above,scale=0.7]{$(0)$};
\filldraw[black] (3.3,-.7) circle (1.3pt)node[above,scale=0.7]{$(2)$};
\draw (0,0) -- (1,0);
\draw (0,0) -- (1,-.7);
\draw (1,-.7) -- (1.5,-.7)node[right]{...};\draw (2,-.7) -- (3.3,-.7);
\draw (0,0) -- (1,.7);
\filldraw[white] (0,.9) circle (.1pt) node[above,black]{$K_3$};
\end{scope}

\begin{scope}[shift={(6.5,0)}]
\filldraw[black] (0,0) circle (1.3pt)node[above,scale=0.7]{$(0)\ $};
\filldraw[black] (1,0) circle (1.3pt)node[above,scale=0.7]{$(0)$};
\filldraw[black] (1,.7) circle (1.3pt)node[above,scale=0.7]{$(0)$};
\filldraw[black] (1,-.7) circle (1.3pt)node[above,scale=0.7]{$(0)$};
\filldraw[black] (2.5,-.7) circle (1.3pt)node[above,scale=0.7]{$(0)$};
\draw (0,0) -- (1,0);
\draw (0,0) -- (1,-.7);
\draw (1,-.7) -- (1.5,-.7)node[right]{...};\draw (2,-.7) -- (2.5,-.7);
\draw (0,0) -- (1,.7);
\filldraw[white] (0,.9) circle (.1pt) node[above,black]{$K_4$};
\end{scope}
\end{tikzpicture}\\
Therefore:
$$d(-M_p,\mathbf t_i)=\max\left\{\frac{c_1(\mathbf s)^2+|\Gamma'|}{4};\mathbf s\in \spinc(W_{\Gamma'}),\mathbf s|_{-M_p}=\mathbf t_i\right\}=
\left\{\begin{array}{cl}0&i=1,2\\ \frac{p-2}{4}&i=3\\ \frac{p+2}{4}&i=4\end{array}\right.
$$

Applying Theorem \ref{thm}, the above computations already give that for distinct $\rot_1,\rot_2$ all structures are tight, and belong to two different isotopy classes, while for equal $\rot_1,\rot_2$ the only tight triples are $(\pm1,\pm1,\mp(p-1))$ and they are isotopic to each other. This proves Claim \ref{example}.

\begin{claim}
Theorem \ref{thm} holds for $M_p$.
\end{claim}

We show this following the two-step analysis described in Section \ref{Sec2}.

\subsection{Detect tightness}
The condition we use to recognize tight structures among all $(M_p,\xi)$ presented by surgery diagrams of Figure \ref{fig:SFS} is an existence of the characteristic covector $c$ as in Theorem \ref{thm:c}.

We give $c$ as $\PD(c)=\alpha h+\sum\alpha_i e_i$ where $\alpha,\alpha_i\in\{\pm1\}$, and such that $(c_i)=\langle c,x_i\rangle=\rot_i-1$. Concretely, the $c$-evaluations on $\Gamma$ belong to one of the following.

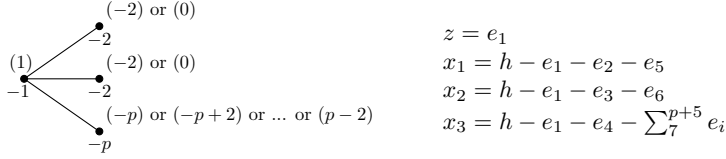
\begin{tikzpicture}
\begin{scope}[shift={(-3.3,0)}]
\filldraw[black] (0,0) circle (1.3pt)node[above,scale=0.7]{$(1)\ $} node[below,scale=0.7]{$-1\ \ $};
\filldraw[black] (1,0) circle (1.3pt)node[above right,scale=0.7]{$(-2)$ or $(0)$} node[below,scale=0.7]{$-2$};
\filldraw[black] (1,.7) circle (1.3pt)node[above right,scale=0.7]{$(-2)$ or $(0)$} node[below,scale=0.7]{$-2$};
\filldraw[black] (1,-.7) circle (1.3pt)node[above right,scale=0.7]{$(-p)$ or $(-p+2)$ or ... or $(p-2)$} node[below,scale=0.7]{$-p$};
\draw (0,0) -- (1,0);
\draw (0,0) -- (1,-.7);
\draw (0,0) -- (1,.7);
\end{scope}
\begin{scope}[shift={(2,0)}]
\filldraw[white] (0,0) circle (.1pt) node[right,black,scale=.9]{$%
\begin{array}{l}z=e_1\\ x_1=h-e_1-e_2-e_5\\ x_2=h-e_1-e_3-e_6\\ x_3=h-e_1-e_4-\sum_7^{p+5}e_i\end{array}$};
\end{scope}
\end{tikzpicture}

Then, for each such $(\alpha,\alpha_i)$ we compute $c|_{\Gamma'}$, and check how its full path ends.

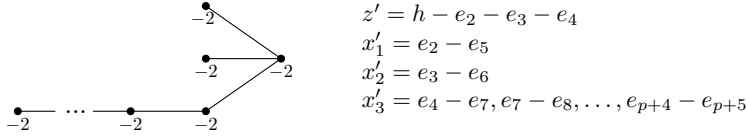
\begin{tikzpicture}
\begin{scope}[shift={(-3.3,0)}]
\filldraw[black] (.5,-.7) circle (1.3pt)node[below,scale=0.7]{$-2$};
\filldraw[black] (2,-.7) circle (1.3pt)node[below,scale=0.7]{$-2$};
\filldraw[black] (3,-.7) circle (1.3pt)node[below,scale=0.7]{$-2$};
\filldraw[black] (3,0) circle (1.3pt)node[below,scale=0.7]{$-2$};
\filldraw[black] (3,.7) circle (1.3pt)node[below,scale=0.7]{$-2\ $};
\filldraw[black] (4,0) circle (1.3pt)node[below,scale=0.7]{$-2$};
\draw (3,-.7) -- (4,0);
\draw (3,0) -- (4,0);
\draw (3,.7) -- (4,0);
\draw (1.5,-.7) -- (3,-.7);
\draw (.5,-.7) -- (1,-.7)node[right]{...};
\end{scope}
\begin{scope}[shift={(1.5,0)}]
\filldraw[white] (0,0) circle (.1pt) node[right,black,scale=.9]{$%
\begin{array}{l}z'=h-e_2-e_3-e_4\\ x_1'=e_2-e_5\\ x_2'=e_3-e_6\\ x_3'=e_4-e_7, e_7-e_8,\dots, e_{p+4}-e_{p+5}\end{array}$};
\end{scope}
\end{tikzpicture}

Below we table all possible $(\alpha,\alpha_i)$ for each given triple $(c_1,c_2,c_3)$. We will make explicit how some $c|_{\Gamma'}$ drop out, and connect the others to the right initial and terminal vector. Also, we will emphasize the appearance of the same characteristic covectors $c|_{\Gamma'}$ in some pairs of $c$-triples.

First observe that (on the level of paths) the order of signs on generators of each leg is unimportant, as they can be shuffled using $\pm2\PD(v')$-steps for $\langle c, v'\rangle=\pm2$. Then there are essentially only two different sign-vectors $(\alpha,\alpha_i)$ for a chosen $c$-triple, differing in the sign of $h$. The two are connected by $\pm2\PD(z')$, applied when $\langle c,z'\rangle=\pm2$. Notice that all these different sign configurations belong to the same surgery presentation.

In the light of the previous paragraph, we record only the number of positive and negative signs on exceptional generators of each leg. Write $\{m+,n-\}_i$ when there are $m$ positive and $n$ negative generators of $L_i$ (counted without $h$ and $e_1$); not to be confused with vectors of signs which record exact sign configuration on corresponding generators. In addition, let $(h+)_{(c_1,c_2,c_3)}$ and $(h-)_{(c_1,c_2,c_3)}$ denote any of sign configurations which belongs to $(c_1,c_2,c_3)$ and has positive, negative respectively, sign on $h$. Look separately at the cases with the same, and later with the distinct $(c_1,c_2)$. 

\begin{small}
$$
\begin{array}{rllllll} c_3\ \ &\underline{(-2,-2,c_3)}&&\ &\underline{(0,0,c_3)}&\\
& (h-) & (h+) && (h-) & (h+)\\
& \{1+,1-\}_1 & \{0+,2-\}_1 && \{2+,0-\}_1 & \{1+,1-\}_1\\
& \{1+,1-\}_2 & \{0+,2-\}_2 && \{2+,0-\}_2 & \{1+,1-\}_2\\
p-2 &\{p+,0-\}_3&\{(p-1)+,1-\}_3&&\{p+,0-\}_3&\{(p-1)+,1-\}_3\\
p-4 &\{(p-1)+,1-\}_3&\{(p-2)+,2-\}_3&&\{(p-1)+,1-\}_3&\{(p-2)+,2-\}_3\\
\vdots&\vdots&\vdots&&\vdots&\vdots\\
-p &\{1+,(p-1)-\}_3&\{0+,p-\}_3&&\{1+,(p-1)-\}_3&\{0+,p-\}_3\\
\end{array}
$$
\end{small}

\noindent
For $(-2,-2,c_3), c_3\in\{p-4,...,-p\},$ there exists a configuration $(h,e_2,e_3,e_4)=(-,-,-,-)$ which drops out: $\langle c,z'\rangle=\langle -h-e_2-e_3-e_4,h-e_2-e_3-e_4\rangle=-4$. Similarly, $(0,0,c_3), c_3\in\{p-2,...,-p+2\},$ drops out at $(h,e_2,e_3,e_4)=(+,+,+,+)$. Therefore, the paths possibly end only for the triples $(-2,-2,p-2)$ and $(0,0,-p)$.

Furthermore, we observe that $(-2,-2,p-2)$ and $(0,0,-p)$ belong to the same full path because configurations $(h-)_{(-2,-2,p-2)}$ and $(h+)_{(0,0,-p)}$ give the same characteristic vector (0's on the third leg, and $(h,e_4):(-,+)\leftrightarrow(+,-)$ with the same evaluation on $z'=h-e_2-e_3-e_4$). This proves also that their (common) path indeed ends, namely at $K_3$ (given by $(h,e_4,e_7,...,e_{p+4},e_{p+5})=(-,-,-,...,-,+)$ for $(0,0,-p)$) on the initial side and at $-K_3$ (as $(h,e_4,e_7,...,e_{p+4},e_{p+5})=(+,+,+,...,+,-)$ for $(-2,-2,p-2)$) on the terminal.

\begin{small}
$$
\begin{array}{rllllll} c_3\ \ &\underline{(-2,0,c_3)}&&\ &\underline{(0,-2,c_3)}&\\
& (h-) & (h+) && (h-) & (h+)\\
& \{1+,1-\}_1 & \{0+,2-\}_1 && \{2+,0-\}_1 & \{1+,1-\}_1\\
& \{2+,0-\}_2 & \{1+,1-\}_2 && \{1+,1-\}_2 & \{0+,2-\}_2\\
p-2 &\{p+,0-\}_3&\{(p-1)+,1-\}_3&&\{p+,0-\}_3&\{(p-1)+,1-\}_3\\
p-4 &\{(p-1)+,1-\}_3&\{(p-2)+,2-\}_3&&\{(p-1)+,1-\}_3&\{(p-2)+,2-\}_3\\
\vdots&\vdots&\vdots&&\vdots&\vdots\\
-p &\{1+,(p-1)-\}_3&\{0+,p-\}_3&&\{1+,(p-1)-\}_3&\{0+,p-\}_3\\
\end{array}
$$
\end{small}

\noindent
Sign configurations adapted to any $c$-triple with distinct $c_1$ and $c_2$ build a connected part of (one of the two) full paths. Indeed, let us see how these parts patch together into a path. For $k\in\{1,...,p-1\}$, we have
$$
(h+)_{(-2,0,p-2k)}\ \underset{c|_{\Gamma'}}{=}\ (h-)_{(0,-2,p-2k-2)}\  \underset{c|_{\Gamma}}{\equiv} (h+)_{(0,-2,p-2k-2)}\ \underset{c|_{\Gamma'}}{=}\ (h-)_{(-2,0,p-2k-4)}
$$
where the first and the last equality denote the same characteristic covector on $\Gamma'$, while the middle equivalence means (different sign distributions of) the same presentation. This separates all characteristic vectors arising from presentations with mixed $(c_1,c_2)$ into two full paths. One starting at 
$K_1$ (as $-h-e_2+e_5+e_3+e_6+e_4+e_7+\cdots+e_{p+5}$ for $(-2,0,p-2)$) and ending at 
$-K_{2}$ (as $+h- e_2-e_5+ e_3-e_6-e_4-e_7-\cdots-e_{p+5}$ for $(-2,0,-p)$) or 
$-K_1$ (as $+h+ e_2-e_5- e_3-e_6-e_4-e_7-\cdots-e_{p+5}$ for $(0,-2,-p)$). 
The other starting at $K_2$ (as $-h+e_2+e_5-e_3+e_6+e_4+e_7+\cdots+e_{p+5}$ for $(0,-2,p-2)$) and ending at 
$-K_{1}$ (as $+h+ e_{2}-e_{5}- e_3-e_6-e_4-e_7-\cdots-e_{p+5}$ for $(0,-2,-p)$) or 
$-K_2$ (as $+h- e_{2}-e_{5}+ e_3-e_6-e_4-e_7-\cdots-e_{p+5}$ for $(-2,0,-p)$). 
The two terminal possibilities depend on parity of $p$ (odd or even). 

In conclusion, translated back into rotation numbers we have obtained the following paths of tight structures, each sharing the same invariants:
\begin{itemize}[leftmargin=.3cm]
	\item $(-1,-1,p-1)$ and $(1,1,-p+1)$ ($\spinc=\mathbf t_4+\mu_a,\ d_3=\frac{2-p}{4}$)
	\item $(-1,1,p-1)$ and $(1,-1,p-3)$ and $(-1,1,p-5)$ and ... ($\spinc=\mathbf t_4+\mu_b\ d_3=0$)
	\item $(1,-1,p-1)$ and $(-1,1,p-3)$ and $(1,-1,p-5)$ and ... ($\spinc=\mathbf t_4+\mu_c,\ d_3=0$)
\end{itemize}

\subsection{Prove overtwistedness and describe contact isotopies}
In our (simplest possible) cases with boundary slopes $\frac{1}{k}, k\in\ZZ$, there is a single continued fraction block for each special fiber. Contact surgery presents direct translation between positive and negative stabilizations (down- and up-cusps) of core Legendrian unknots and positive and negative basic slices in the decomposition of a continued fraction block with slopes $-1$ and $-k$. The generators forming the corresponding leg (and by that, the dual vertices) in the plumbings above can be thought of as another way of layering solid torus into $k$ slices.

We need contact topological interpretation for the steps in full paths.

First, the unimportance of sign permutations in the legs coincide with the shuffling of basic slices within a single continued fraction block \cite[Lemma 4.14]{Ho.I}. Moreover, \cite[Section 6]{GLS} provides sufficient isotopy moves between contact structures presented by different surgery diagrams. Let us spell this out. Since the moves in \cite{GLS} are given by the matrices of signs whose coefficients are $q_j^i$, the number of positive basic slices in the $j^\text{th}$ continued fraction block of the $i^\text{th}$ leg, in our case only $(q_0^1, q_0^2, q_0^3)$, we rewrite previously obtained paths of tight structures in this language, changing rotation numbers to $q_0^i$'s:
\begin{itemize}
	\item $(0,0,p-1)$ and $(1,1,0)$
	\item $(0,1,p-1)$ and $(1,0,p-2)$ and $(0,1,p-3)$ and ...
	\item $(1,0,p-1)$ and $(0,1,p-2)$ and $(1,0,p-3)$ and ...
\end{itemize}

Now, we notice that conditions which caused a full path to drop out, and so prevented our tightness criterion to work, exactly agree with condition for which overtwistedness can be proved. And finally, there are contact isotopies between pairs of surgery presentations which share the same path. Let us recite.

\begin{prop}\cite[Propositions 6.3, 6.1 \& 6.4]{GLS}
Let a contact structure on $M_p$ be given by $(q_0^1, q_0^2, q_0^3)$ as above. Then the triples $(1,1,q_0^3)$ with $q_0^3\neq 0$ and $(0,0,q_0^3)$ with $q_0^3\neq p-1$ present overtwisted structures. Between other presentations, there are following contact isotopies:
$$
(1,0,q_0^3)\simeq \left\{\begin{array}{l} (0,1,q_0^3+1) \mbox{ when }q_0^3<p-1\\ (0,1,q_0^3-1) \mbox{ when }q_0^3>0\end{array}\right.\ \ \mbox{and}\ \ (1,1,0)\simeq(0,0,p-1).\ \ \ \ \ \square
$$
\end{prop}

\subsection*{Problems in general}
Examples shown above are special in several ways. In general, it can happen that the full path associated to some presentation $(c_v)_{v\in\Gamma}$ drops out, although all characteristic covectors computed from $(\alpha,\alpha_i)$-configurations which restrict to $(c_v)_{v\in\Gamma}$ satisfy the bounds $v\cdot v\leq\langle c,v\rangle\leq-v\cdot v\ \text{ for all } v\in\Gamma'$. Also, not all the steps in a full path need to be presentable, that is, arising from some tuple of rotation numbers.  (For  examples of such paths, look at the two ``applications" in \cite{LS}.) That said, we need to find out how the (subsequent) presentations of the same path are related, when neither of their characteristic covectors on $\Gamma'$ coincides (Corollary \ref{cor:paths}). Finally, we need new conditions for overtwistedness (Proposition \ref{prop:ot}) and isotopies (Proposition \ref{prop:isotopies}), which will explain such behavior of full paths.

\section{Characteristic covectors, tightness, and full paths}\label{Sec4}

In Subsection \ref{Step1}, we have associated characteristic covectors on $\Gamma'$ to any given surgery presentation. Here we investigate their full paths. Namely, how these paths end, and which presentations share the same path. In order to do so, we first observe that certain $2\PD$-steps do not change the presentation (Subsection \ref{Ss1}). Then, we explore the only remaining central step -- concretely, we explain it on the level of homology generators (Subsection \ref{Ss2}). In the following Subsections, we are then concerned with the associated change in $c|_\Gamma$, whether this new $c|_\Gamma$ comes from some presentation and when it leads to the end of the path (Subsection \ref{Ss3}). Moreover, we describe (in Subsection \ref{Ss4}) the first presentable $c|_\Gamma$ (or the end of path) following any possible starting point.

\begin{notation}
We describe a characteristic $2$-cohomology class $c\in H^2(R;\ZZ)$ as $\PD(c)=\alpha h+\sum\alpha_ie_i$ where $\alpha,\alpha_i\in\{\pm1\}$. In the following, vectors of signs correspond to parts of the coefficient-vector $(\alpha,\alpha_i)$, covering generators of (usually) a single $\Gamma$- or $\Gamma'$-vertex. 

To a single vertex we often refer by its self-intersection. When a vertex is written out in generating classes, these are called starting, middle and last, according to the position; explicitly, if $v=e_s-\sum_{j=s+1}^le_j$, then $e_s$ is starting, $e_l$ is last, and all others are middle. On legs, the starting generator of a vertex and the last generator of the previous vertex coincide.

Presentability will be assigned to dual vectors and it means that the corresponding manifold side arises from a contact presentation, that is, the manifold-side evaluations can be expressed by rotation numbers as in Figure \ref{fig:graphs}.
\end{notation}

Let our starting point be a characteristic vector which comes from a contact presentation, and which satisfies $v\cdot v\leq\langle c,v\rangle\leq-v\cdot v$ (otherwise we have already dropped out). We will follow the path only in one direction -- towards the initial vector. Recall that the corresponding step is given by $-2\PD(v)$ for some $v$ with $\langle c,v\rangle=v\cdot v$, and the vector we aim at satisfies $v\cdot v+2\leq\langle c,v\rangle\leq-v\cdot v$. Everything could be verbatim repeated with opposite signs in the direction of the terminal vector.

\subsection{Steps on legs}\label{Ss1}

First, we observe that steps taken for $v\neq z'$ never change the presentation considered, neither the path drops out at any of these vertices. (To the remaining case $v=z'$ we dedicate Subsection \ref{central}.)

\begin{lem}
Characteristic vectors $c$ and $c- 2\PD(v), v\neq z', \langle c,v\rangle= v\cdot v,$ always belong to the same surgery presentation.
\end{lem}

\proof
As these vertices ($v\in\Gamma', v\neq z'$) are described by $v=e_i-\sum e_j$, the evaluation of characteristic covector $c$ reaches the self-intersection when presenting generators all admit the same sign as in the vertex. So, $-2\PD(v)$ changes their signs from $(+-\cdots-)$ to $(-+\cdots+)$. But this change has no effect on the evaluation of $c$ on any of the $\Gamma$-vertices.

Indeed, from the way how the exceptional classes are chosen we see that each $e_j$ starts some new vertex, either one on the manifold side or one on the dual side. So, the starting and the last generator of $v$ are non-starting on the manifold side, while all its middle generators are starting (and last) generators of manifold vertices. Hence, the restriction of $c$ to the generators of $v$ evaluates trivially on $\Gamma$, $\langle c|_v,\Gamma\rangle=0$, and is therefore independent of sign.

Since these (manifold-side) evaluations directly correspond to rotation numbers, by neither of these moves we switch between presentations.
\endproof

\begin{lem}
All drop-outs occur in the center $z'=h-e_2-e_3-e_4$ of the dual star.
\end{lem}

\proof
We notice that all the vertices in legs of $\Gamma'$ are formed by exactly as many generators ($e_j$'s) as the value of their self-intersections. Hence, there is no way to drop out at any of them. So, the only possible drop-out happens at $z'$ when the signs of generators $h$ and $e_2,e_3,e_4$ are all the same ($(++++)$ or $(----)$, and $\langle c,z'\rangle=\pm4$).
\endproof

In sum, we may assume the initial condition $v\cdot v+2\leq\langle c,v\rangle\leq-v\cdot v$ is violated only at the central vertex $z'$ -- such vector can be easily reached by finishing all possible $-2\PD$-steps on legs, which either sweep out the problem or transfer it to the center. (As each $-2\PD$-step pushes the problem to the neighboring vertices, we are successively completing the steps, as long as we do not run into a vertex $v$ which despite of the $-2$-change does not evaluate as $v\cdot v$, or we reach the end of the leg.) In particular, neither non-central vertex is of the form $(+-\cdots -)$.

\subsection{Central step}\label{central}\label{Ss2}

After the above reduction, covector $c$ either drops out at $z'$, presents the initial vector or, it reaches self-intersection at $z'$. For the latter, the generators forming $z'=h-e_2-e_3-e_4$ take values: either $(+---)$ or $(-+--)$, up to reordering the legs. The $-2\PD$-step taken next changes exactly these generators by twice $(-+++)$. In the first case we stay in the same presentation, as the step only switches the signs in the pairs $(h,e_i), i=2,3,4$, preserving the evaluation on all the influenced manifold vertices, $x_1,x_2,x_3$. In the second case, we can (on the level of generators) instead of simply adding $-2z'$ to the given description of $\PD(c)$, first change the sign configuration, without changing the dual $c|_{\Gamma'}$ and with controlled (seen later) change on the manifold side $c|_\Gamma$, and then do the $-2\PD$-step as above, not influencing the manifold side.

\begin{alg}[Central step or turn]\label{alg}
Whenever we arrive, after possibly renumbering the legs, at $c$ with $(h,e_2,e_3,e_4)=(-+--)$ and $\langle c,v\rangle\neq v\cdot v$ for all $v\neq z'$, the next step in the full path is given by the characteristic covector $\overline{c}$ as follows. Denoting vertices of $L_1'$ by $\{v_0',\dots,v_{k_1'}'\}$ and their generators as $v_i'=e_1^i-\sum_{j=2}^{l_i}e_j^i$ with $e_{l_i}^i=e_1^{i+1}$ and $e_1^0=e_2$, define $\PD(\overline{c})=\PD(c)+2h-2e_2$ and modify it as follows:
$$\begin{array}{lll}
    \text{for }i\in\{0,\dots,k_1'\}  & \text{if }\langle \overline{c},e_1^i\rangle\neq\langle c,e_1^i\rangle: & \text{for }j\in\{2,\dots,l_i\} \\
    &&\ \ \ \ \ \text{if }\langle c,e_j^i\rangle=+1: \\
    &&\ \ \ \ \ \PD(\overline{c})=\PD(c)-2e_j^i\ \&\text{ stop} \\
    & \text{if }\langle \overline{c},e_1^i\rangle=\langle c,e_1^i\rangle: & \text{stop.}
\end{array}$$
Then add $-2z'$ to so obtained sign configuration $\PD(\overline{c})$.
\end{alg}

To prove well-definedness, we need that such reformulation always exists (the inner loop in our Algorithm always stops, Lemma \ref{lem:exist}) and that uniqueness, ensured by always taking the first positive generator (chosen ordering of the inner loop), can be explained by the independence of order, at least as far as contact presentations are concerned (Lemma \ref{lem:unique}). 

\begin{lem}\label{lem:exist}
Every characteristic vector $c_{\Gamma'}$ with $\PD(c|_{z'})=(-+--)$ can be achieved by another distribution of signs, with positive sign on $h$; it is associated to a different manifold vector (possibly non-presentable).
\end{lem}

\proof
Starting at the center $z'$, the two distributions are given by $(-+--)$ and $(+---)$. The switch of the $h$-sign with the opposite sign of $e_2$, does not impose any change into the second and the third dual leg. For the first leg, the appropriate adaptation of signs, which results in the same dual evaluation, exists because of exclusion of any $(+-\cdots-)$-configurations (that is, the assumption $\langle c,v\rangle\neq v\cdot v\ \text{ for any } v\neq z'$).
\endproof

\begin{lem}\label{lem:unique}
As a sign on one middle generator of a dual vertex is changed, all of them need to be changed (independent of order) before we get back into presentable. A turn of the last generator can result in a presentable vector only when all prior middle generators are negative.
\end{lem}

\proof 
For a  covector to be presentable, all dual vertices have to have same-signed middle generators, because these generators on the manifold side are forming a chain of $-2$'s, zero being their only possible rotation number.

For the second claim, suppose on the contrary the middle signs on some $v'$ are positive. Changing the sign of its last generator (from positive to negative) forces a switch of all the signs in the following chain (if any) of dual $-2$'s (to preserve dual evaluations). Then, this influences the evaluation on the next non-$(-2)$ dual vertex $w'$, which can be corrected by changing one of its later generators from positive to negative. If the middle generators of $w'$ are already negative or if we get them all negative by the current turn, we have found (independent of further changes) a manifold-side vertex which starts at positive (second last) generator in $v'$ and has all further signs negative. If by the change of one middle generator not all of them are negative, the vector is non-presentable by the first part. If all generators of $w'$ are positive, and we turn the last one, we need to repeat the same argument with $w'$ in place of $v'$. It remains to check whether we could get presentable result by correcting only starting and last generators of all following (necessarily, fully positive) dual vertices. But if not before, the process ends in non-presentable, giving $(+-\cdots -)$ on the last manifold vertex.
\endproof 

To sum up, the central turns are the only significant steps in following possible changes on manifold vectors, and by that, in presentations.  We may assume that after each central turn also all $-2\PD$-steps on legs are finished. 

\subsection{On turning sequences and presentability}\label{Ss3}

To begin, notice how to recognize the ends of a full path.

\begin{lem}\label{lem:1-2-3}
If after a central step, covector $c$ on the starting dual vertices evaluates as their self-intersection, $\langle c,x_i'\rangle=x_i'\cdot x_i'$:
\begin{itemize}[leftmargin=.7cm]
  \item on at most one leg, we have arrived at the initial end;
  \item on two legs, the full path continues;
  \item on all three legs, this causes a drop-out.
\end{itemize}
\end{lem}

\proof
The maximal starting dual evaluations tell us on how many legs we need further $-2\PD$-steps. The evaluation on $z'$, $\langle c, z'\rangle$, right after a central turn is $+2$. If further turns are needed for one leg only we do not reach $-2$ central evaluation again and the corresponding vector is initial; with two we get back to $\langle c, z'\rangle=-2$ and we continue with another central turn; three gives a drop-out in $(h,e_2,e_3,e_4)=(----)$.
\endproof

\begin{cor}\label{cor:full+}
A presentation $\xi$, whose $\mathcal{P}_\xi$ properly ends, necessarily admits a leg, starting in a fully positive vertex. (If presentation corresponds to the initial vector, there are two fully positive starting vertices.)
\end{cor}

\proof
For $\PD(c)$ take a sign configuration which evaluates on manifold vertices according to the rotation numbers of $\xi$, which takes minus sign on $h$, and for which $\langle c,v'\rangle\neq v'\cdot v'\ \text{ for all } v'\in\Gamma'\backslash\{x_i';i=1,2,3\}$. (This is the stage right after a central turn.) As in Lemma \ref{lem:1-2-3} above, there is a leg, say $L_1$, for which $\langle c,x_1'\rangle\neq x_1'\cdot x_1'$. We prove that on this leg $\langle c,x_1\rangle=a_0^1-2$ holds, that is, the generators of $x_1$ (apart from $h,e_1$) are positive.

Write out $x_1$ as $h-e_1-e_2-e_5-\cdots-e_{J}$. The signs on the generators up to $e_{J-1}$ are positive as otherwise we would have shuffled the negative sign to $e_5$ by $-2\PD$-steps on consecutive dual vertices of square $-2$ (resulting in $\langle c, x_1'\rangle=-2$ for $x_1'=e_2-e_5$). The positivity of $e_J$ follows from presentability via the following claim.

\begin{claim*} A presentable covector on neither dual vertex  takes the form $(+-\cdots-+)$.\end{claim*}
\proof A proof of this fact is basically the same as the second part of the proof of Lemma \ref{lem:unique}. Suppose on the contrary, there is such a dual vertex; it is not the last vertex of the dual leg, because it would give the last vertex on the manifold side with self-intersection $-2$ and $+2$ $c$-evaluation. But then, every non-$(-2)$ dual vertex further on the dual leg needs to have again negative middle signs (otherwise we have found a manifold vertex, starting in the negative sign of the previous non-$(-2)$ with all following generators positive) and positive last one (because of $(+-\cdots-)$ exclusion). After all, ending in the impossible last dual vertex.\endproof

Since also $(+-\cdots-)$-configuration on any dual vertex, except $x_i',$ is excluded by $\langle c,v'\rangle\neq v'\cdot v'$, and since middle generators of any dual vertex are same-signed, we get that $e_J$ is positive. It is a middle generator of a dual vertex starting in positive $e_{J-1}$.
\endproof

The leg with fully positive starting vertex is the one which in the reordering of Algorithm \ref{alg} takes role of $L_1$. When we wish to emphasize according to which leg in the actual structure the central step is done, we refer to it as a turn of $L_i$.

Since the evaluation of characteristic covector on $L_i$-vertices changes only by turns of $L_i$, we may separately study their influence.

\begin{lem}\label{lem:oneleg}
Let $c$ be a presentable non-initial characteristic covector. Assume that it evaluates on the vertices of some leg $L=(-a_0,-a_1,\dots,-a_j,-a_{j+1},\dots,-a_k)$ as follows: $$\langle c,L\rangle=(a_0-2, a_1-2,\dots,a_j-2,a_{j+1}-2-2n_{j+1},\dots,a_k-2-2n_k)$$
where $k\geq j, n_{j+1},\dots,n_{k}\geq0$ and $n_{j+1}>0$.

The path runs into the next possibly presentable covector $\bar{c}$ only after $$1+1+(a_1-1)+(a_2-1)(a_1-1)+\cdots+(a_{j-1}-1)\cdots(a_1-1) \text{ turns of } L,$$ in: $$\langle \bar{c},L\rangle=(-a_0, -a_1+2,\dots,-a_j+2,a_{j+1}-2n_{j+1},\dots,a_k-2-2n_k).$$
\end{lem}

\proof
To be illustrative, we explicitly write out all the generators involved in the first few turns. Below are the two sides, $L_1$ and $L_1'$, in homology generators; the $\ast$-symbol stands for truncation only.\\

\noindent{\small $\begin{array}{llll}
 L_1:\qquad  x_1= h-e_1-e_2-e_5-\cdots-e_{J-1}- & e_J &               &                     \\
                                                                       & e_J-&e_{J+1}   &                     \\
                                                                       &  	      &e_{J+1}- & e_{J+2}         \\
                                                                       &	      &		      &  \ddots           \\
                                                                       &        &               & e_{K-1}-e_{K}-e_{K+1}-\ast
\end{array}$\\
$\begin{array}{llll}                                                              
 L_1':\qquad    x_1'=  e_2- & e_5         & & \\
                                 & \ddots      & & \\
                                 & e_{J-2}-   & e_{J-1} & \\
                                 &                 & e_{J-1}-e_J-e_{J+1}- \cdots -&e_{K} \\
                                 &&&\ast
\end{array}$}

In this notation, the starting part of $L_1$ and the evaluation of $c$ on it take values:
$$L_1=(-J+3,\underbrace{-2,\dots,-2}_{K-J-1},-T,-S,\ast) \text{ and } \langle c,L_1\rangle=(J-5,0,\dots,0,M,N,\ast).$$

By the first turn, according to the Algorithm, we change generators up to $e_{J}$ -- it does not influence further dual vertices, but a new vector can be presentable only when all the middle generators $e_J,...,e_{K-1}$ are same-signed. Therefore, in order to (possibly) reach presentable vector again we have to repeat turning of this particular leg $(K-J)$-times. Resulting manifold vector is of the form $(-J+3,0\dots,0,M+2,N,\ast)$, its presentability depends on the $(+2)$-changed manifold vertex $e_{K-1}-e_K-e_{K+1}-\ast$.

In terms of generators, we have reached another presentation exactly when $e_K$-sign is negative. The positive $e_K$-sign, on the other hand, requires another turn, but this forces some further changes to preserve the dual. Namely, we need to change signs on generators of the following chain of $-2$'s, and one (without loss of generality, first) middle generator afterwards. The resulting vector is not necessarily presentable, provided the starting point was, it depends on presentability of the vertex starting in the (last changed) middle generator ($+2$ rotation change). But if it is, the new presentation is $(-J+3,0,\dots,0,-T+2,N+2,\ast)$; for this, we need to turn this leg $(K-J+1)$-times.

Continuing in the same manner, we trace similar behavior at all levels. Concretely. We are successively turning fully positive vertices, which influences the evaluation on the following manifold vertex by $+2$. If the result is presentable, we have finished. Otherwise, the following vertex was also fully positive, at the moment its evaluation is minus self-intersection, and it will have turned under the influence of another turn of the previous vertex. For that we need to bring the previous vertex back to maximal rotation, using (again) influence of the previous vertices on leg. But notice that each vertex is influenced only by turns of the vertex just before it. Therefore, to come from maximal rotation through minus self-intersection to minimal rotation on some vertex $v_{k+1}$, we need to influence it by two turns of its first previous vertex $v_k$. This in turn is obtained by $(a_k-1)$ turns of its previous vertex $v_{k-1}$, by first to get from minus self-intersection to minimal rotation, and then by the step of $+2$ to maximal rotation. This explains the number of steps and finishes the proof.
\endproof

Obviously, the leg (its vertices with self-intersections) together with the sign configuration (in presentable, rotation numbers) determine when the leg is turned. In particular, it specifies the gaps between the subsequent turnings of the same leg, when some other leg needs to be turned in order for the path to continue. Actually, the reverse also holds.

\begin{lem}\label{lem:periods}
A form of a leg together with a distribution of signs on its generators is completely described by the sequence of its turns. 
\end{lem}

\proof As before, separately state (and argue for) the first step.

\begin{claim*}
Between two subsequent turnings of the same leg $L$ there are always either $a_0-2$ or $a_0-1$ turnings of other two legs. \end{claim*}

\proof
Remember that all generators (but possibly last) of the starting manifold vertex on the turning leg are positive. Since by each turn of other legs we change starting evaluation by $+2$ (through the change of $h$-sign from negative to positive), the gap is determined by the number of generators of the starting vertex. Its variation by one is due to whether the dual vertex following $-2$'s is also fully negative after the $L$-turn.
\endproof

That said, given a turning sequence, we get $a_0$ out of the size of gap between subsequent turnings. If the gap is always the same ($a_0-1$), this means that $v_0$ is the only vertex on $L$, and its self-intersection is $-a_0$.

Call this number of central turns between two turns of $L$, the ($0^\text{th}$) period of $L$. Let us define higher periods for a turning sequence:
\begin{itemize}[leftmargin=1.9cm]
  \item[$1^\text{st}$ period:] number of times the $0^\text{th}$ period is $a_0-2$ before it turns $a_0-1$;
  \item[$2^\text{nd}$ period:] number of times the $1^\text{st}$ period is $a_0'-1$ before it turns $a_0'-2$;
  \item[$k^\text{th}$ period:] number of times the $(k-1)^\text{th}$ period is $a_{k-2}'-1$ before it turns $a_{k-2}'-2$.
\end{itemize}
The numbers should not be read from the first time round. As suggested by notation, they correspond to the self-intersections of dual vertices, hence they determine $L'$, and by that $L$. The initial distribution of signs can be now recognized from the values of periods before the first change.
\endproof

\subsection{Restrictions on the whole structure}\label{Ss4}

We look at all possible (presentable) entries. For each we continue its path as long as it reaches another presentation, or otherwise it ends, either by a drop-out or a (non-presentable) initial vector. Throughout we assume that the (normalized) Seifert constants are ordered $r_1\geq r_2\geq r_3$. In order to reduce possibilities we invoke the $L$-space condition.

\subsubsection*{$L$-space condition}
Recall the numerical condition for $M$ to be an $L$-space: there are no coprime integers $m,a$ such that $\frac{1}{r_1}>\frac{m}{a}, \frac{1}{r_2}>\frac{m}{m-a}, \frac{1}{r_3}>m$; we say that coefficients $(r_1,r_2,r_3)$ are not realizable. As a direct consequence of this condition we observe that:
\begin{enumerate}[leftmargin=.7cm, label=(\roman*)]
  \item\label{i} $r_1\geq \frac{1}{2}$, equivalently one leg starts with $-2$ (otherwise the realizability condition is satisfied for coprime $m=2,a=1$), and more
  \item\label{ii} if $L_1=(\underbrace{-2,\dots,-2}_{k},\ast)$ (then $\frac{1}{r_1}>\frac{k+2}{k+1}=-(\underbrace{[2,\dots,2]}_{k+1})$), then $\frac{1}{r_2}\leq k+2$, equivalently $x_2\cdot x_2\geq -k-2$.
\end{enumerate}

\subsubsection*{Dual configurations} In the following arguments, there will frequently appear a pair of (truncated) legs which are dual to each other, that is, describing a lens space and its dual. Recall that the coefficients of the two are related as follows (here, $-2^{\times b_\iota}$ means a chain of $b_\iota$-many $-2$'s): $$\begin{array}{llccccc}L_i&=&(-b_1-2,&-2^{\times b_2},& -b_3-3,&\dots,&-b_m-2)\\ L_j&=&(-2^{\times b_1},&-b_2-3,&-2^{\times b_3},&\dots,& -2^{\times b_m})\end{array}.$$ The inverses of the continued fractions they describe, add up to $-1$.

\medskip
So, our starting point is a presentable characteristic covector, which does not present an initial end, or a drop-out. Thus, exactly one leg of the corresponding presentation starts in a fully positive vertex. We separate the cases: in Proposition \ref{prop:1&2} we gather presentations for which either $v_0^1$ or $v_0^2$ is stabilized fully positively, and in Proposition \ref{prop:3} we cover presentations for which $v_0^3$ is fully positively stabilized.

\begin{prop}\label{prop:1&2}
Let $c$ be a presentable non-initial characteristic covector, associated to a presentation with fully positive starting vertex on $L_i$, either $L_1$ or $L_2$. This means that it evaluates on the vertices of  $$L_i=(-a_0^i,-a_1^i,\dots,-a_j^i,-a_{j+1}^i,\dots,-a_{k_i}^i)$$ as in Lemma \ref{lem:oneleg}: $$\langle c,L_i\rangle=(a_0^i-2, a_1^i-2,\dots,a_j^i-2,a_{j+1}^i-2-2n_{j+1}^i,\dots,a_{k_i}^i-2-2n_{k_i}^i)$$
for some $0\leq j\leq k_i, n_{j+1}^i,\dots,n_{k_i}^i\geq0$ and $n_{j+1}^i>0$.

Denote coefficients on other two legs by $a_\kappa^\lambda$, and let $c$ evaluate as $\langle c,v_0^3\rangle=a_0^3-2-2n_0^3$ on the first vertex of the third leg $L_3$, and $$\langle c,L_l\rangle=(-a_0^l, -a_1^l+2,\dots,-a_k^l+2,a_{k+1}^l-2-2n_{k+1}^l,\dots,a_{k_l}^l-2-2n_{k_l}^l)$$ on the leg $L_l, l\neq i,3$, for some $-1\leq k\leq k_l,\ n_{k+1}^l\leq a_{k+1}^l-3$.

Furthermore, define $m$ and $N_{m+1}^l$ as follows:
$$\begin{array}{l}m:=\max\{\kappa\leq k;\ \exists N: \text{ denominator of }[a_0^l,\dots,a_\kappa^l,N]\leq n_0^3\}\\
N_{m+1}^l:=\left\{\!\!\!\begin{array}{l}n_{k+1}^l+1,\text{ if denominator of }[a_0^l,\dots,a_k^l,n_{k+1}^l+1]^{-1}\leq n_0^3, \text{ or }\\ N\in [1,a_{m+1}^l): \!\!\!\!\!\begin{array}{l}\text{ denominator of }[a_0^l,\dots,a_m^l,N]^{-1}\leq n_0^3 \text{ and }\\ \text{ denominator of }[a_0^l,\dots,a_m^l,N+1]^{-1}>n_0^3, \text{ otherwise.}\end{array}\end{array}\right.\end{array}$$

Then the full path of $c$ behaves as follows.
\begin{enumerate}[leftmargin=.6cm]
\item If $-[a_0^i,\dots,a_{j-1}^i]^{-1}-[a_0^l,\dots,a_m^l,N_{m+1}^l]^{-1}<1,$ the full path drops out.

\item If $-[a_0^i,\dots,a_{j-1}^i]^{-1}-[a_0^l,\dots,a_m^l,N_{m+1}^l]^{-1}=1 \text{ and } m=k, N_{m+1}^l=n_{k+1}^l+1$, we reach a new presentation $\bar{c}$ which on the three legs takes the following values:
$$\begin{array}{ccl}\langle \bar{c},L_i\rangle&=&(-a_0^i, -a_1^i+2,\dots,-a_j^i+2,a_{j+1}^i-2n_{j+1}^i,\dots,a_{k_i}^i-2-2n_{k_i}^i)\\
\langle \bar{c},L_l\rangle&=&(a_0^l-2, a_1^l-2,\dots,a_{k+1}^l-2,a_{k+2}^l-2n_{k+2}^l,\dots,a_{k_l}^l-2-2n_{k_l}^l)\\
\langle \bar{c},v_0^3\rangle&=&a_0^3-2-2n_0^3+2 D\end{array}$$
where $D$ stands for the denominator of $[a_0^l,\dots,a_k^l,n_{k+1}^l+1]^{-1}$, and the evaluations on the rest of $L_3$ remain the same as for $c$.

\item Otherwise, the path continues in non-presentable and reaches the non-presentable initial end.
\end{enumerate}
\end{prop}

\begin{remark}
Rewrite the coefficients up to $a_j^i$ in the $b$-notation used for dual configurations above, so for appropriate $b_\iota\geq0$ (notice $b_1>0$ on $L_1$):
$$L_i=(-2^{\times b_1},-b_2-3,-2^{\times b_3},\dots, -2^{\times b_J},-a_{j+1}^i,\dots,-a_{k_i}^i).$$
Also for $L_l$, truncated as continued fractions in the Proposition, take $$(a_0^l,\dots,a_m^l,N_{m+1}^l)=(-b_1'-2,-2^{\times b_2'},-b_3'-2,\dots,-b_J'-2).$$ Then the conditions, given in the Proposition in terms of continued fraction sums, can be restated as:
\begin{itemize}[leftmargin=0.7cm]
\item the two continued fractions add up to $1$ when $b_k'=b_k\ \text{ for all } k$;
\item the sum is greater than $1$ when for $K=\min\{k;b_k'\neq b_k\}:\ b_K'<b_K$ if $K$ odd, $b_K'>b_K$ if $K$ even;
\item the sum is smaller than $1$ when for $K=\min\{k;b_k'\neq b_k\}:\ b_K'>b_K$ if $K$ odd, $b_K'<b_K$ if $K$ even.
\end{itemize} 
\end{remark}

\proof[Proof of Proposition \ref{prop:1&2}]
We need to observe how specific behavior of the path restricts possible forms of a covector, and by that, of a presentation.

To meet another presentation, recall that we need to swap the signs of all generators forming the fully positive vertices $v_0^i,\dots,v_j^i$ (Lemma \ref{lem:oneleg}). To achieve this, we need certain number of $L_i$-turns, which are arranged in the turning sequence, uniquely determined by the form of $L_i$. So, any turn of other two legs should appear at exactly specified non-turning stages of $L_i$. 

Now notice that any turn of $L_3$ (before finishing the specified sequence) would immediately end the path (in non-presentable). Indeed, it has to appear after $k$ (or $k+1$) turns of $L_1$ (which is, the $1^\text{st}$ period of $L_1$). So after it $L_2$ needs to be turned (its starting coefficient being $a_0^2\leq k+2$, hence its $0^\text{th}$ period being at most $k+1$) and also $L_1$ needs to be turned (having $0^\text{th}$ period $0$ or $1$). But this already means we have arrived at the initial end, see Lemma \ref{lem:1-2-3}.

Therefore, the turning sequence of $L_i$ (up to its $(j+1)^\text{th}$ vertex) exactly specifies the turning sequence of $L_l$ (as far that presentability on $L_i$ is reached), its turns being in non-turning points of $L_i$, and vice versa. It is exactly the turning sequence of the dual leg with all evaluations fully negative. Rewritten in terms of relations between continued fractions, the two legs are of dual forms if and only if the corresponding continued fractions add up to one. In the dual (negative) leg for the last entry the significant information is the number of negative signs, as turning sequence depends only on whether we have reached maximal evaluation.

Taken together, we have obtained.
\begin{enumerate}[leftmargin=0.7cm]
\item If not all gaps in the turning sequence of $L_i$ are filled by turns of $L_l$, and at the same time, the sequence is not quit by the turn of $L_3$ before or at the time when first such non-filled gap appears, then the full path drops out.

\item If the two continued fractions add up to exactly one, this means the turning sequences of corresponding legs exactly fit together, and we reach another presentation -- if and only if the evaluation on the starting vertex of the third leg is negative enough, not to quit the sequence of turnings interchanging between $L_i$ and $L_l$. That is, there has to be more negative generators as there are turns of $L_i$ and $L_l$, which equals the denominator of corresponding continued fractions. Form of the new presentation is determined by Lemma \ref{lem:oneleg}.

\item Otherwise means that, either we hit into some turning point of $L_l$ before reaching the next gap in the sequence of $L_i$, or we reach a turning point of $L_3$ at or before the time when $L_i$-gap is not filled by $L_l$-turn for the first time. As observed above, in these cases, the full path properly ends with the initial vector, but it is necessarily non-presentable because we have not yet reached the first possibly presentable stage as specified in Lemma \ref{lem:oneleg}.
\endproof\end{enumerate}

\begin{prop}\label{prop:3}
Throughout the path, there can be at most two non-initial characteristic covectors for which the starting vertex of $L_3$ is fully positive, that is $\langle c,v_0^3\rangle=a_0^3-2$.

If $\langle c,v_1^3\rangle\neq a_1^3-2$ or if $L_3=(v_0^3)$, the turn of $L_3$ is presentable, the two presentations differ in: $\langle \bar{c},v_0^3\rangle=-a_0^3,\ \langle \bar{c},v_1^3\rangle=\langle c,v_1^3\rangle+2,$ and for $l=1,2,\ \langle \bar{c},v_0^l\rangle=\langle c,v_0^l\rangle+2$.

If $\langle c,v_1^3\rangle= a_1^3-2$, and $c$ is not terminal, the turn of $L_3$ necessarily makes the continuation of the path non-presentable, and ends it in a non-presentable initial end.

If $\langle c,v_1^3\rangle= a_1^3-2$, and $c$ is terminal, let us write out the $c$-evaluations at the terminal end:
$$\begin{array}{l}  \langle c,L_1\rangle=(-a_0^1, -a_1^1+2,\dots,-a_j^1+2,a_{j+1}^1-2-2n_{j+1}^1,\dots,a_{k_1}^1-2-2n_{k_1}^1)\\ \langle c,L_2\rangle=(-a_0^2, -a_1^2+2,\dots,-a_k^2+2,a_{k+1}^2-2-2n_{k+1}^2,\dots,a_{k_2}^2-2-2n_{k_2}^2) \end{array}$$
for some $0\leq j\leq k_1$ and $0\leq k\leq k_2$. Then:
\begin{enumerate}[leftmargin=.6cm]
\item If for maximal $J\leq j, K\leq k$ such that $-[a_0^1,\dots,a_{J}^1,3]^{-1}-[a_0^2,\dots,a_K^2,2]^{-1}=1,$ the denominator of the two fractions is smaller than $a_0^3$, then the full path drops out.

\item If $\langle c,v_2^3\rangle\neq a_2^3-2$, and there exist $J\leq j, K\leq k$ such that $$\begin{array}{ll}-[a_0^1,\dots,a_{J}^1,3]^{-1}-[a_0^2,\dots,a_K^2,2]^{-1}=1 \text{ and } & n_{J+1}^1\geq 2\text{ or }L^1=(v_0^1,\dots,v_J^1) \\ & n_{K+1}^2\geq 1\text{ or }L^2=(v_0^2,\dots,v_K^2) \end{array},$$ and the denominator of the two fractions equals $a_0^3$, we reach a new presentation $\bar{c}$ which on the three legs takes the following values:
$$\begin{array}{ccl}
\langle \bar{c},L_1\rangle&=&(a_0^1-2, a_1^2-2,\dots,a_J^1-2,\langle c, v_{J+1}^1\rangle+4, \langle c,v_\iota^1\rangle_{\iota=J+2}^{k_1})\\
\langle \bar{c},L_2\rangle&=&(a_0^2-2, a_1^2-2,\dots,a_{K}^2-2,\langle c, v_{K+1}^2\rangle+2, \langle c,v_\kappa^2\rangle_{\kappa=K+2}^{k_2})\\
\langle \bar{c},L_3\rangle&=&(-a_0^3,-a_1^3,\langle c, v_2^3\rangle+2, \langle c,v_\mu^3\rangle_{\mu=3}^{k_3}).\end{array}
$$
This $\bar{c}$ presents the initial end of the full path.

\item Otherwise, the path continues in non-presentable and reaches the non-presentable initial end.
\end{enumerate}
\end{prop}

\proof
As observed in the proof of Proposition \ref{prop:1&2}, any time when path runs into a presentation with fully positive $v_0^3$ (not at its terminal end), it reaches the initial vector, either before or after a turn of $L_3$. Therefore, if the characteristic vector before the $L_3$-turn is non-initial and presentable, the only other presentation, which can appear as we continue the path, can occur straight after this turn. The resulting vector is presentable if and only if $\langle c,v_1^3\rangle\neq a_1^3-2$ (when exists). Relation between the two presentations is as always read from Lemma \ref{lem:oneleg} (the simplest possible -- one-turn -- case).

The only remaining option is to have a presentable terminal end with $\langle c,v_0^3\rangle=a_0^3-2$. In that case the turn of $L_3$ does not end the path, and if this turn is not presentable itself, we need to look for any possible following presentation. Since presentability on $L_3$ can be recovered only by a turn of $L_3$, and since according to above this turn ends the path, we might meet such a presentation (only) at the initial end. This in particular means that turns of $L_1$ and $L_2$ in between the two turns of $L_3$ should begin and end with a vector which is presentable on these two legs. The first turn after the $L_3$-turn, and the last turn before another $L_3$-turn are done according to $L_1$. 

As before, we inductively determine that, in order for the turning sequences of $L_1$ and $L_2$ to fit together (being interchangeably turned until the second $L_3$-turn), the two legs need to have fully negative starting vertices, forming almost dual vectors. Almost in the sense that there is no ``last pair", that is, the two vectors as given in the paragraph on dual configurations end by $-b_m-3$ and $-2^{\times b_m}$ instead of $-b_m-2$ and $-2^{\times b_m}$. In other words, they are dual when enlarged by $3$ and $2$, respectively.

The three possibilities are now given as before:
\begin{enumerate}[leftmargin=.7cm]
\item We have not reached the turning point of $L_3$ yet, but the sequence of turnings of $L_1$ and $L_2$ cannot continue.
\item The turn of $L_3$ appears exactly in the moment when neither $L_1$ nor $L_2$ can be turned, and the sign configuration on them is presentable. Additionally, we need that $\langle c,v_2^3\rangle\neq a_2^3-2$ to reach presentability of $L_3$ as well.
\item Otherwise, either the first two legs hit into a common turning point, the third leg finishes the sequence early not having enough negative generators, or (simply) the terminal vector obtained as in (2) is not presentable because $\langle c,v_2^3\rangle= a_2^3-2$ (hence, after the $L_3$-turn, $\langle \bar{c}, v_2^3\rangle=a_2^3$). In all the cases, the path stops in non-presentable initial end.
\endproof\end{enumerate}

Above we described how the successive presentations in the path are related to each other and indicate what property causes a drop-out. Any given presentation can be now either walked through these stages to the proper ends of the full path or it drops out. Joining results (taking into account also their analogues obtained by following the path in the terminal direction) we obtain the following picture. Here, the presentations are given as evaluations of characteristic covectors on generators of $H_2(W_\Gamma)$, written as triples of vectors $c^i$ whose entries are $c_j^i=\langle c,v_j^i\rangle$. Vectors are truncated -- we write out only the relevant part and hide the rest into $\ast$.

\begin{cor}[Full path components]\label{cor:paths}
If a given presentation $\xi$ does not admit both a fully positive and a fully negative starting vertex, its full path drops out at $\xi$. Moreover, a full path drops out when it runs into a presentation given by either of the following characteristic covectors $c|_{\Gamma}$, independently of how the three vectors continue in the hidden $\ast$-part.\\
For some $(i,l)\in\{(1,2),(2,1)\}$: 
$$c|_\Gamma=\left(\begin{array}{c} a_0^i-2\\ a_1^i-2\\ \vdots\\ a_j^i-2\\ a_{j+1}^i-2-2n_{j+1}^i\\ \ast\end{array}
\left|\begin{array}{c} -a_0^l\\ -a_1^l+2\\ \vdots\\ a_{k+1}^l-2-2n_{k+1}^l\\ a_{k+2}^l-2-2n_{k+2}^l\\ \ast\end{array}\right|
\begin{array}{c} a_0^3-2-2n_0^3\\ \ast\end{array}\right)$$
for which $-[a_0^i,\dots,a_{j-1}^i]^{-1}-[a_0^l,\dots,a_m^l,N_{m+1}^l]^{-1}<1$ holds for $N_{m+1}^l$ defined as in Proposition \ref{prop:1&2}.\\
Or:
$$c|_\Gamma=\left(\begin{array}{c} -a_0^1\\ -a_1^1+2\\ \vdots\\ -a_j^1+2\\ a_{j+1}^1-2-2n_{j+1}^1\\ \ast\end{array}
\left|\begin{array}{c} -a_0^2\\ -a_1^2+2\\ \vdots\\ -a_k^2+2\\ a_{k+1}^2-2-2n_{k+1}^2\\ \ast\end{array}\right|
\begin{array}{c} a_0^3-2\\ a_1^3-2\\ \ast \end{array}\right)$$
such that for maximal $J\leq j, K\leq k$ with $-[a_0^1,\dots,a_{J}^1,3]^{-1}-[a_0^2,\dots,a_K^2,2]^{-1}=1,$ the denominator of the two fractions is smaller than $a_0^3$.

In the terminal direction, symmetrically, a drop-out occurs at presentations with oppositely stabilized surgery link (that is, surgery diagrams given by the same but reversely oriented link).

\medskip
Any two presentations $\xi_1,\xi_2$ whose associated characteristic vectors meet at the same path $\mathcal{P}_{\xi_1}=\mathcal{P}_{\xi_2}$ are related by the sequence of rotation number changes, each taking one of the following forms. The pairs are presented in form of $c|_\Gamma$ and they have to be identical on all further generators, hidden in $\ast$. \\
Either for $(i,l)\in\{(1,2),(2,1)\}$:
$$\begin{array}{l}
\left(\begin{array}{c}a_0^i-2\\ a_1^i-2\\ \vdots\\ a_j^i-2\\ a_{j+1}^i-2-2n_{j+1}^i\\ \ast\end{array}
\left|\begin{array}{c}-a_0^l\\ -a_1^l+2\\ \vdots\\ a_{k+1}^l-2-2n_{k+1}^l\\ a_{k+2}^l-2-2n_{k+2}^l\\ \ast\end{array}\right|
\begin{array}{c}a_0^3-2-2n_0^3\\ \ast \end{array}\right)\ \simeq\ \\ \\
\simeq\ \left(\begin{array}{c}-a_0^i\\ -a_1^i+2\\ \vdots\\ -a_j^i+2\\ a_{j+1}^i-2n_{j+1}^i\\ \ast\end{array}
\left|\begin{array}{c}a_0^l-2\\ a_1^l-2\\ \vdots\\ a_{k+1}^l-2\\ a_{k+2}^l-2n_{k+2}^l\\ \ast\end{array}\right|
\begin{array}{c}a_0^3-2-2(n_0^3-D)\\ \ast\end{array}\right)
\end{array}$$
where $D$ is the denominator of $[a_0^l,\dots,a_{k}^l,n_{k+1}^l+1]^{-1}$, and $k, n_{k+1}^l$ satisfy $$1=-[a_0^i,\dots,a_{j-1}^i]^{-1}-[a_0^l,\dots,a_k^l,n_{k+1}^l+1]^{-1}.$$

\noindent Or:
$$\begin{array}{l}
\left(\begin{array}{c}a_0^1-2-2n_0^1\\ \ast\\ \end{array} \left|\begin{array}{c}a_0^2-2-2n_0^2\\ \ast\\ \end{array}\right| \begin{array}{c}a_0^3-2\\a_1^3-2-2n_1^3\\ \ast\end{array}\right)\ \simeq\\ \\ 
\simeq\ \left(\begin{array}{c}a_0^1-2n_0^1\\ \ast\\ \end{array} \left|\begin{array}{c} a_0^2-2n_0^2\\ \ast\\ \end{array}\right| \begin{array}{c}-a_0^3\\a_1^3-2n_1^3\\ \ast\end{array}\right)\end{array}.$$

\noindent Or:
$$\begin{array}{l}
\left(\begin{array}{c}  -a_0^1\\ -a_1^1+2\\ \vdots\\ -a_J^1+2\\ \vdots\\ -a_j^1+2\\ a_{j+1}^1-2-2n_{j+1}^1\\ \ast\end{array}  
\left|\begin{array}{c}-a_0^2\\ -a_1^2+2\\ \vdots\\ -a_K^2+2\\ \vdots\\ -a_k^2+2\\ a_{k+1}^2-2-2n_{k+1}^2\\ \ast\end{array}\right| \begin{array}{c} a_0^3-2\\ a_1^3-2\\ a_2^3-2-2n_2^3\\ \ast \end{array}\right)\ \simeq\\ \\
\simeq\ \left(\begin{array}{c}  a_0^1-2\\ \vdots \\a_J^1-2\\a_{J+1}^1+2-2n_{J+1}^1\\ \ast\end{array}  
\left|\begin{array}{c}a_0^2-2\\ \vdots\\ a_K^2-2\\ a_{K+1}^2-2n_{K+1}^2\\ \ast\end{array}\right|  
\begin{array}{c}-a_0^3\\ -a_1^3+2\\ a_2^3-2n_2^3 \\ \ast \end{array}\right)
\end{array}$$
where for $J\leq j, K\leq k,$  $$\begin{array}{ll}-[a_0^1,\dots,a_{J}^1,3]^{-1}-[a_0^2,\dots,a_K^2,2]^{-1}=1 \text{ and } & n_{J+1}^1\geq 2\text{ or }L^1=(v_0^1,\dots,v_J^1) \\ & n_{K+1}^2\geq 1\text{ or }L^2=(v_0^2,\dots,v_K^2) \end{array},$$ and the denominator of the two fractions equals $a_0^3$.\hfill $\square$
\end{cor}

\section{Convex surface theory, overtwistedness, and isotopies}\label{Sec5}

To prove that the isotopic classification of tight structures is contained in  the full paths of their dual covectors, we need to observe that presentations sharing the same path are indeed isotopic, and relate drop-outs to overtwistedness. This section covers part of the proof outlined in Subsection \ref{Step2}.

To begin with, remember two simple properties of full paths, which have direct convex theoretic interpretation. The first is the shuffling property of basic slices within a single continued fraction block \cite[Subsubsection 4.4.5]{Ho.I}, which can be in Heegaard Floer interpretation recovered by $2\PD$ steps on the consecutive dual vertices of square $-2$. The second is, necessary condition for tightness, that the presentation contains both a leg starting in a fully positive vertex, and a leg starting in a fully negative vertex. In full paths, a fully positive starting vertex is required by Corollary \ref{cor:full+}, and a fully negative one by its analogue when following the path in terminal direction. In convex surface theory, other presentations can be seen to fail the conditions of the Gluing Lemma \cite[Theorem 4.25]{Ho.I}, as in \cite[Proposition 6.3]{GLS}, but can be also understood as a special case of overtwistedness proved below.

Let us now state the result as predicted from the Heegaard Floer picture, as in Corollary \ref{cor:paths}. We encode contact presentations into ``matrices of negative signs", that is, triples of vectors $q^i$, possibly of different length, whose coefficients are $q_j^i$, the number of negative basic slices in the $j^{\text{th}}$ continued fraction block of the $i^{\text{th}}$ singular fiber. The three vectors in Propositions are truncated, so that we write out only the relevant part (on which overtwistedness is decided, or which behaves non-trivially under isotopy moves) and hide the rest into $\ast$. Analogously, we can define ``matrices of positive signs". Notice that the ones counting negative slices directly correspond to the relations obtained in the initial direction of the full path in Section \ref{Sec4}. With positive slices they correspond to symmetric relations in terminal direction. To describe isotopy moves it is enough to give only pairs of matrices of negative (or only positive) signs, while to encode conditions for overtwistedness, the two are different.

\begin{prop}[Overtwistedness conditions]\label{prop:ot}
Let a contact presentation be described by either of the following matrices of signs, negative or positive:
\begin{enumerate}[leftmargin=0.7cm, label=O\arabic*.]
\item\label{O1} For some $(i,l)\in\{(1,2),(2,1)\}$: 
$$(q^i | q^l | q^3)=\left(\begin{array}{c}0\\ 0\\ \vdots\\ 0 \\ n_{j+1}^i \\ \ast \end{array}
\left|\begin{array}{c}a_0^l-1\\ a_1^l-2\\ \vdots\\ a_k^l-2\\ n_{k+1}^l\\ n_{k+2}^l\\ \ast\end{array}\right|
\begin{array}{c}n_0^3\\ \ast \end{array}\right)$$
for which $-[a_0^i,\dots,a_{j-1}^i]^{-1}-[a_0^l,\dots,a_m^l,N_{m+1}^l]^{-1}<1$ holds for $N_{m+1}^l$ defined as in Proposition \ref{prop:1&2}.
\item\label{O2} $$(q^1 | q^2 | q^3)=\left(\begin{array}{c} a_0^1-1\\ a_1^1-2\\ \vdots\\ a_j^1-2\\ n_{j+1}^1\\ \ast\end{array} 
\left|\begin{array}{c} a_0^2-1\\ a_1^2-2\\ \vdots\\ a_k^2-2\\ n_{k+1}^2\\ \ast\end{array}\right|
\begin{array}{c} 0\\ 0\\ \ast \end{array}\right)$$
such that for maximal $J\leq j, K\leq k$ with $-[a_0^1,\dots,a_{J}^1,3]^{-1}-[a_0^2,\dots,a_K^2,2]^{-1}=1,$ the denominator of the two fractions is smaller than $a_0^3$.
\end{enumerate}
Then, independently of the basic slice decompositions of further continued fraction blocks ($\ast$-part of vectors), the corresponding contact structure is overtwisted.
\end{prop}

\begin{prop}[Isotopy conditions]\label{prop:isotopies}
The following pairs of matrices give isotopic contact structures, provided all coefficients are in the range $n_j^i\in[0,a_j^i-2], n_0^i\in[0,a_0^i-1]$, and the further basic slice decompositions ($\ast$-parts) are the same.
\begin{enumerate}[leftmargin=0.7cm, label=I\arabic*.]
\item\label{I1} On $(q^i | q^l | q^3)$ for $(i,l)\in\{(1,2),(2,1)\}$:
$$
\left(\begin{array}{c} 0\\ 0\\ \vdots\\ 0\\ n_{j+1}^i\\ \ast\end{array}
\left|\begin{array}{c} a_0^l-1\\ a_1^l-2\\ \vdots\\ a_k^l-2\\ n_{k+1}^l\\ n_{k+2}^l\\ \ast\end{array}\right|
\begin{array}{c} n_0^3\\ \ast\end{array}\right)\ \simeq\ 
\left(\begin{array}{c} a_0^i-1\\ a_1^i-2\\ \vdots\\ a_j^i-2\\ n_{j+1}^i-1\\ \ast\end{array}
\left|\begin{array}{c} 0\\ 0\\ \vdots\\ 0\\ 0\\ n_{k+2}^l-1\\ \ast\end{array}\right|
\begin{array}{c} n_0^3-D\\ \ast\end{array}\right)
$$
where $D$ is the denominator of $[a_0^l,\dots,a_{k}^l,n_{k+1}^l+1]^{-1}$, and $k, n_{k+1}^l$ satisfy $$1=-[a_0^i,\dots,a_{j-1}^i]^{-1}-[a_0^l,\dots,a_k^l,n_{k+1}^l+1]^{-1}.$$

\item\label{I2} On $(q^1 | q^2 | q^3)$: $$\left(\begin{array}{c} n_0^1\\ \ast\\ \end{array} \left|\begin{array}{c} n_0^2\\ \ast\\ \end{array}\right| \begin{array}{c} 0\\ n_1^3\\ \ast\end{array}\right)\simeq
\left(\begin{array}{c}n_0^1-1\\ \ast\\ \end{array} \left|\begin{array}{c} n_0^2-1\\ \ast\\ \end{array}\right| \begin{array}{c}a_0^3-1\\ n_1^3-1\\ \ast\end{array}\right).$$

\item\label{I3} On $(q^1 | q^2 | q^3)$: $$
\left(\begin{array}{c} a_0^1-1\\ a_1^1-2\\ \vdots\\ a_J^1-2\\ \vdots\\ a_j^1-2\\ n_{j+1}^1\\ \ast\end{array}
\left|\begin{array}{c} a_0^2-1\\ a_1^2-2\\ \vdots\\ a_K^2-2\\ \vdots\\ a_k^2-2\\ n_{k+1}^2\\ \ast\end{array}\right|
\begin{array}{c} 0 \\ 0\\ n_2^3 \\ \ast \end{array}\right)\ \simeq\ 
\left(\begin{array}{c} 0\\ \vdots \\ 0\\ n_{J+1}^1-2\\ \ast\end{array}
\left|\begin{array}{c} 0\\ \vdots \\ 0\\ n_{K+1}^2-1\\ \ast\end{array}\right|
\begin{array}{c} a_0^3-1\\ a_1^3-2\\ n_2^3-1 \\ \ast \end{array}\right)
$$
where for $J\leq j, K\leq k,$  $$\begin{array}{ll}-[a_0^1,\dots,a_{J}^1,3]^{-1}-[a_0^2,\dots,a_K^2,2]^{-1}=1 \text{ and } & n_{J+1}^1\geq 2\text{ or }L^1=(v_0^1,\dots,v_J^1) \\ & n_{K+1}^2\geq 1\text{ or }L^2=(v_0^2,\dots,v_K^2) \end{array},$$ and the denominator of the two fractions equals $a_0^3$.
\end{enumerate}
\end{prop}

Proof of both Propositions is postponed till the end of the section, after a note on contact topological foundations of isotopies and some general computation of slopes. 

\subsection{State traversals and contact isotopies}

Convex decomposition of Seifert fibration we are working with consists of the neighborhoods of the three singular fibers $F_i$ and the background circle bundle over the pair of pants. To ensure the product structure in the complement of $F_i$'s we use non-normalized coefficients $M(0;r_1-1,r_2,r_3)$. The results here rely on Honda's classification of tight structures on separated pieces, namely solid tori \cite{Ho.I} and circle bundles over the pair of pants $\Sigma$ \cite[Subsection 5.1]{Ho.II}. In the context of isotopies, we are mostly concerned with the changes of boundary slopes obtained by thickening tubular neighborhoods of singular fibers $F_i$. The corresponding picture is the $3$-punctured sphere with given boundary slopes.

Concretely, we recognize correspondences between presentations relying on the following lemma.

\begin{lem} \cite[Lemma 4.13]{GSch} \label{lem:A}
Let $\Sigma$ be a pair of pants and $\xi$ a tight contact structure on $\Sigma\times S^1$ whose boundary $-\partial(\Sigma\times S^1)=T_1\cup T_2\cup T_3$ consists of tori in standard form with $\#\Gamma_{T_i}=2$ for $i=1,2,3$, and slopes $s(T_1)=-\frac{p_1}{q}, s(T_2)=-\frac{p_2}{q}, s(T_3)=\infty$. Suppose that there exists a pair of pants $\Sigma'\subset\Sigma$ such that $\Sigma\times S^1$ decomposes as $\Sigma\times S^1=\Sigma'\times S^1\cup C_1\cup C_2,\ C_i=T_i\times I$, with $\xi|_{C_i}$ minimally twisting and where $\xi|_{\Sigma'\times S^1}$ is a tight contact structure with infinite boundary slopes such that the section $\Sigma'\times\{\theta\}$ for some $\theta\in S^1$ is convex with dividing set consisting of arcs, each connecting two different boundary components.

If $s(T_2)=-\frac{p_2}{q}<0$ and both $\xi|_{C_i}$ decompose into basic slices of the same sign, there exists a convex annulus $A$ bounded by the Legendrian rulings of $T_1$ and $T_2$, and without boundary parallel dividing curves. \hfill $\square$
\end{lem}

In our case, the decomposition $\Sigma\times S^1=\Sigma'\times S^1\cup C_1\cup C_2$ always exists as we are dealing with the zero-twisting structures; both thickened tori are minimally twisting \cite[Lemma 5.1]{Ho.II}. The condition on the background structure, no boundary parallel dividing curves on the convex section, follows (as in \cite[Lemma 5.4]{GLS}) from the fact that $\xi$ is appropriate on $\Sigma\times S^1$ (such that, there is no embedded $T^2\times I$ with $T^2$ isotopic to a boundary component, and $I$-twisting at least $\pi$). And the latter is satisfied for any $\Sigma\times S^1$, cut as a background out of tight small Seifert manifold \cite[Lemma 2.4]{Wu}. With addition of \cite[Section 3]{G}, Lemma \ref{lem:A} can be reformulated in the sense of \cite[Lemmas 5.7, 5.8]{GLS}.

\begin{lem}\label{lem:L1L2}
Let $\Sigma$ be a pair of pants and let $\xi$ be an appropriate contact structure on $\Sigma\times S^1$ with convex boundary $-\partial(\Sigma\times S^1)=T_1\cup T_2\cup T_3$ with $\#\Gamma_{T_i}=2$ for $i=1,2,3$, and boundary slopes $s(T_1)=-\frac{p_1}{q}, s(T_2)=-\frac{p_2}{q}, s(T_3)=\infty$.
\begin{enumerate}[leftmargin=.8cm, label=(L\arabic*)]
  \item\label{L1} If there exists a collar neighborhood $C_3\subset\Sigma\times S^1$ of $T_3$, which is minimally twisting with boundary slopes $\infty$ and $\frac{p_1+p_2-1}{q}$, whose basic slices are all same-signed, and for which $\xi|_{(\Sigma\times S^1)-C_3}$ coincides with the unique tight structure with boundary slopes $-\frac{p_1}{q},-\frac{p_2}{q},\frac{p_1+p_2-1}{q}$, and maximal twisting $-q$, then signs of basic slices in the decomposition of $C_{1}$ and $C_2$ are all opposite to $C_3$-signs.
  \item\label{L2} And conversely, if $C_1$ and $C_2$ decompose into same-signed basic slices, then there exists $C_3$ composed of opposite-signed slices, with boundary slopes $\infty$ and $\frac{p_1+p_2-1}{q}$, and such that its complement is a unique tight structure as above. 
\end{enumerate}
\end{lem}

\proof[Proof (following \cite{GLS})]
For \ref{L1}, uniqueness of a tight structure with given properties is stated in \cite[Proposition 3.3]{G}. The fact that the signs in the decomposition of collars are opposite, can be read from the relative Euler class evaluation on vertical annuli $A_i\subset C_i$, which have boundaries in vertical Legendrian divides on $\infty$-side and in Legendrian rulings on the other boundary. If we complete these annuli with annuli through the pair of pants up to $T_3$ for $A_1, A_2$, and up to $T_1, T_2$ for $A_3$, we get two pairs of homologically equivalent, but oppositely oriented, annuli. As the Euler class evaluation on all the extended parts is zero (first having boundary in Legendrian divides, second living inside $-q$-maximal twisting), the evaluation on $A_1$ and $A_2$ is opposite to that of $A_3$. Therefore, the collars $C_1, C_2$ decompose into basic slices all of the same sign, opposite to the signs in $C_3$.

For \ref{L2}, we take the unique tight $\Sigma''\times S^1$ (as described in \ref{L1}) and attach to it a thickened torus with slopes $\frac{p_1+p_2-1}{q},\ \infty$, and slices signed oppositely to the ones in $C_1$. Then according to \ref{L1} the signs on collars in decomposition $\Sigma'\times S^1\cup C_1\cup C_2$ are again opposite. And we have built up a contact structure, isotopic to original in all three pieces ($\Sigma'\times S^1$ has same dividing set on the pair of pants, while $C_1$ and $C_2$ have the same Euler class evaluations).
\endproof

\subsection{Slicing and continued fractions}

We give a short reflection on the slopes of glued-up torus and its slicing. Denote $V_i$ the standard convex neighborhood of $F_i$ with boundary $-\partial(M\backslash V_i)$ trivialized by $1\choose 0$ the horizontal direction of $\Sigma\times 1$ and $0\choose 1$ the direction of fiber, and from the other side $\partial V_i$ by the meridian and some longitude. The last being chosen so that translation $A_i:\partial V_i\rightarrow-\partial(M\backslash V_i)$ is given by $\ A_i=$ {\tiny$\begin{pmatrix}     \alpha_i & \alpha_i'   \\  -\beta_i   & -\beta_i' \end{pmatrix}$} where $\frac{\beta_i}{\alpha_i}=r_i\ (r_1-\ 1 \text{ for the first leg})$; in terms of the continued fraction expansion we have $-\frac{\alpha_i}{\beta_i}=[a_0^i,...,a_{k_i}^i], -\frac{\alpha_i}{\alpha_i'}=[a_{k_i}^i,...,a_0^i], -\frac{\beta_i}{\beta_i'}=[a_{k_i}^i,...,a_1^i]$. Now, the $\infty$-slope of a thickened neighborhood $U_i$ of a singular fiber corresponds to $[a_{k_i}^i,...,a_0^i]$ in the torus basis, and the slopes of the factorization can be obtained in order (from outside in) by decreasing the last entry of this fraction.

We will be interested in slopes of tori, which peel off certain sequences of basic slices from thickened neighborhoods $U_i$, and their expression in the background basis. Notice the following general behavior. 

\begin{lem}\label{lem:slopes}
The slope of torus which peels off $\sum_0^{j-1}(a_\iota^i-2)+m$ outer basic slices from $U_i$, as seen from $-\partial(M\backslash V_i)$, is independent of inner continued fraction blocks in the decomposition of $U_i$, that is, of vertices $a_j^i,\dots,a_{k_i}^i$, farther down the legs. It equals $[a_0^i,\dots,a_{j-1}^i,m]^{-1}$ ($[a_0^i,\dots,a_{j-1}^i,m]^{-1}+1$ in case $i=1$).
\end{lem}

\proof
The slope of interest is in the torus basis expressed as $[a_{k_i}^i,..., a_j^i-m]$.

Recall the matrix form of a negative continued fraction \begin{center} $[a_{k_i}^i,\dots,a_0^i]\leftrightarrow$ {\tiny $\left[\begin{pmatrix}    -a_{k_i}^i & 1   \\  -1   &  0\end{pmatrix}^{-1}\cdots\begin{pmatrix}     -a_{0}^i & 1   \\  -1   &  0\end{pmatrix}^{-1}\right]^2=\left[\begin{pmatrix}     -\beta_i' & -\alpha_i'  \\   \beta_i & \alpha_i \end{pmatrix}\right]^2$},\end{center}
and notice that it is exactly the inverse of our identification $A:\partial V_i\rightarrow-\partial(M\backslash V_i)$.

 Hence,  we get the desired slope in the second column of:
 
{\tiny
$$\begin{pmatrix} -a_0^i & 1   \\  -1   &  0 \end{pmatrix}
\cdots\underbrace{\begin{pmatrix}     -a_{j+1}^i & 1   \\  -1   & 0 \end{pmatrix}
\cdots\begin{pmatrix}     -a_{k_i}^i & 1   \\  -1   & 0 \end{pmatrix}
\begin{pmatrix}    -a_{k_i}^i & 1   \\  -1   & 0 \end{pmatrix}^{-1}
\cdots\begin{pmatrix}     -a_{j+1}^i & 1   \\  -1   &  0\end{pmatrix}^{-1}}_{I}
\begin{pmatrix}     -a_j^i+m & 1   \\  -1   &  0\end{pmatrix}^{-1}=
$$

$$=\begin{pmatrix}    -a_0^i  & 1   \\  -1   & 0 \end{pmatrix}
\cdots\begin{pmatrix}  -a_j^i    & 1   \\  -1   & 0 \end{pmatrix}
\begin{pmatrix}     0 & -1   \\  1   &  -a_j^i+m\end{pmatrix} =
\begin{pmatrix}    -a_0^i  & 1   \\  -1   & 0 \end{pmatrix}
\cdots\begin{pmatrix}  -a_{j-1}^i    & 1   \\  -1   & 0 \end{pmatrix}
\begin{pmatrix}     1 & m   \\  0   & 1 \end{pmatrix} 
$$
}
Indeed, independent of $a_\iota^i, \iota\geq j$.

Now, comparing the second columns:

{\tiny
$$
\left[\begin{pmatrix}    -a_0^i  & 1   \\  -1   & 0 \end{pmatrix}
\cdots\begin{pmatrix}  -a_{j-1}^i    & 1   \\  -1   & 0 \end{pmatrix}
\begin{pmatrix}     1 & m   \\  0   & 1 \end{pmatrix} \right]^2=:\begin{pmatrix}      A   \\  B \end{pmatrix}\leftrightarrow
\begin{pmatrix}     -B   \\  -A \end{pmatrix}=:
\left[\begin{pmatrix}    0  & -1   \\  1   & -a_0^i \end{pmatrix}
\cdots\begin{pmatrix}  0   & -1   \\  1   &   -a_{j-1}^i\end{pmatrix}
\begin{pmatrix}     0 & -1   \\  1   &  -m\end{pmatrix}\right]^2,
$$
}
we express the slope as $[a_0^i,\dots,a_{j-1}^i,m]^{-1}\ ([a_0^i,\dots,a_{j-1}^i,m]^{-1}+1$ for $i=1$).
\endproof

This independence of inner layers, allows us to compute background-basis slope of any sequence of slices (from outside in) on a truncated leg. In the opposite direction, if the slope of peeled-off slices is $[a_0^i,\dots,a_{j-1}^i,m]^{-1}\ ([a_0^i,\dots,a_{j-1}^i,m]^{-1}+1$ if $i=1$), this in torus basis corresponds to $[a_k^i,\dots,a_{j}^i-m]$ when $m\leq a_j^i-1$. When $m=a_j^i$, we get $[a_k^i,\dots,a_{j+1}^i,0]$, undefined as continued fraction, but in terms of the chain of surgeries, the $0$-framed meridian cancels $a_{j+1}$ which results in $[a_k^i,\dots,a_{j+2}^i]$.

\subsection{Proofs}

Proofs of Propositions \ref{prop:ot} and \ref{prop:isotopies} are stated for matrices of negative signs, but they can be verbatim repeated for positive ones. Without loss of generality, basic slices within each continued fraction block are shuffled so that the negative slices are outer.

\proof[Proof of Proposition \ref{prop:ot}]
Guiding principle is as follows. Look at the two singular tori $U_\mu, U_\nu$ whose outermost slices are negative. If we can peel such sequences of negative basic slices from $U_\mu, U_\nu$ that their inner boundary tori $T_\mu, T_\nu$ have slopes with the same denominator, say $-\frac{p_\mu}{q}, -\frac{p_\nu}{q}$, we can use Lemma \ref{lem:A} to find a torus parallel to $\partial U_\lambda$ of slope $\frac{p_\mu+p_\nu-1}{q}$, call it $T$. Whenever this slope is not greater than the critical slope of the singular fiber, $\text{Crit}(F_\lambda)$ (that is, the slope of meridian of the glued-up torus in the background basis), there exists a torus of critical slope between $\partial U_\lambda$ and $T$, which proves overtwistedness (see also \cite[Section 2]{LS.II}). Furthermore, if the slope of $T$ is such that the thickened torus between $\partial U_\lambda$ and $T$ (whose basic slices are all positive by \ref{L2}) forms a basic slice together with some of the slices in the original decomposition of $U_\lambda$, any negative basic slice in this glued-together basic slice implies overtwistedness by the Gluing Lemma. Below we analyze the slopes in each of the cases separately.

\ref{O1} Consider first the structures of the first kind with $(i,l)=(1,2)$. Around $F_2$ there are only negative slices in first $k+1$ continued fraction blocks and $n_{k+1}^2$ of them in the block corresponding to the vertex $v_{k+1}^2$ (here we shuffle them to be its outer). Around $F_3$, on the other hand, we have $n_0^3$ negative slices (again, shuffled so that they are the outer). Thus, peeling off from $U_2$ basic slices up to the $(N_{m+1}^2)^\text{th}$ slice of $-a_{m+1}^2$-block, we obtain the slope $[a_0^2,\dots,a_m^2,N_{m+1}^2]^{-1}$ (in the background basis), which can be joined by cutting annulus to the torus with slope $-\frac{1}{D}$ around $F_3$ (peeling off $D$ slices from $U_3$, where $D$ is denominator of $[a_0^2,\dots,a_m^2,N_{m+1}^2]^{-1}$, which is at most $n_0^3$ by assumption). That way, we have found a torus $T$ parallel to $T_1$ of slope $-[a_0^2,\dots,a_m^2,N_{m+1}^2]^{-1}$. 

Now, observe that the critical slope of $F_1$ is between $$[a_0^1,\dots,a_{j-1}^1-1]^{-1}+1\leq\text{Crit}(F_1)=[a_0^1,\dots,a_{k_1}^1]^{-1}+1\leq[a_0^1,\dots,a_{j-1}^1]^{-1}+1.$$ Our assumed condition gives $-[a_0^2,\dots,a_m^2,N_{m+1}^2]^{-1}<[a_0^1,\dots,a_{j-1}^1]^{-1}+1$. If also a bit more holds true, $-[a_0^2,\dots,a_m^2,N_{m+1}^2]^{-1}\leq\text{Crit}(F_1)$, the torus $T$ embraces the critical one. Otherwise, we have $-[a_0^1,\dots,a_{k_1}^1]^{-1}-[a_0^2,\dots,a_m^2,N_{m+1}^2]^{-1}\geq 1$ and we can truncate both fractions so that the truncations add up to exactly one \cite[Lemma 3.2]{LL}. So, $-[a_0^1,\dots,a_{J}^1]^{-1}-[a_0^2,\dots,b_M^2]^{-1}= 1$ for some $J\in\{j,\dots,k_1\}$ and $M\in\{0,\dots,m+1\}$, and $b_M^2=a_M^2$ for $M\leq m$ or $b_M^2=N_{m+1}^2$ for $M=m+1$. Peeling off from $U_2$ only slices of first $M+1$ outer blocks and corresponding (as many as the denominator of $[a_0^2,\dots,b_M^2]^{-1}$, which is certainly less than or equal to $D\leq n_0^3$) slices in $U_3$, the slope of $T$ is $[a_0^1,\dots,a_{J}^1]^{-1}+1$ in the background basis. By the text under Lemma \ref{lem:slopes} this equals $[a_{k_1}^1,\dots,a_{J+2}^1]$ in torus basis, and $T$ bounds a basic slice with a torus $T_1$ of slope $[a_{k_1}^1,\dots,a_{J+2}^1-1]$. For tightness, the conditions of the Gluing Lemma require for all the subslices of a glued-together basic slice to be positive, but this is not satisfied as the toric annulus bounded by $T_1$ and $T$ contains $({j+1})^\text{th}$ continued fraction block ($J\geq j$) and $n_{j+1}^1>0$ by assumption.

For $(i,l)=(2,1)$, the arguments are the same, but here the induced slope of $T$ (built from peeling-off tori in $U_1$ and $U_3$) equals $[a_0^1,\dots,a_m^1,N_{m+1}^1]^{-1}+1$, while the critical slope is given by $$[a_0^2,\dots,a_{j-1}^2-1]^{-1}\leq\text{Crit}(F_2)= [a_0^2,\dots,a_{k_2}^2]^{-1}\leq[a_0^2,\dots,a_{j-1}^2]^{-1}.$$

\ref{O2} Structures of the second kind admit negative basic slices in the outer layers of $U_1$ and $U_2$. The background-basis slopes $[a_0^1,\dots,a_J^1,3]^{-1}+1$ on $T_1$ around $F_1$ and $[a_0^2,\dots,a_K^2,2]^{-1}$ on $T_2$ around $F_2$ -- which add up to zero -- are reached by peeling off the corresponding sequences of (negative) basic slices when $v_{J+1}^1, v_{K+1}^2$ exist, and by decreasing the twisting number of the Legendrian singular fibers $F_1$ or $F_2$ by stabilizing when $L_1=(v_0^1,\dots,v_J^1)$ or $L_2=(v_0^2,\dots,v_K^2)$. Joining the two tori $T_1$ and $T_2$ by an annulus interpolating between the rulings, and edge-rounding, we obtain a torus $T$ around $F_3$ of slope $-\frac{1}{D}$ where $D$ is the denominator of $[a_0^1,\dots,a_J^1,3]^{-1}$. By assumption, the denominator $D$ is not greater than $a_0^3-1$, hence the obtained slope $-\frac{1}{D}$ is smaller than or equal to $-\frac{1}{a_0^3-1}<\text{Crit}(F_3)$, that is, $T$ embraces the critical torus.
\endproof

\proof[Proof of Proposition \ref{prop:isotopies}]
Recognition of isotopies, in all cases, follows the same steps. First, we apply \ref{L2} to get additional thickened torus $C$ around the singular torus $U_\lambda$ with positive outermost slices -- traverse outer layers, whose basic slices are all negative, from $U_\mu$ and $U_\nu$. This new collar together with some, say $n$, continued fraction blocks around $F_\lambda$ join into a thicker basic slice, separated into positive basic subslices. The isotopy is now given by reversing all these signs. In case we have used all continued fraction blocks of $U_\lambda$, it is interpreted in destabilization followed by opposite stabilization of a core knot. Otherwise, $C$ together with $n+1$ outermost continued fraction blocks in $U_\lambda$ builds a continued fraction block. Its signs can be shuffled, resulting in the $+2$-change in its innermost  $(n+1)^\text{th}$ block (one negative slice replaced by positive) and turn of sign on all basic slices that form $C$ and the first $n$ continued fraction blocks (from positive to negative). The basic slices around the other two fibers, $F_\mu$ and $F_\nu$, are then adapted according to \ref{L1} -- the peeled off ones change their signs from negative to positive, others remain untouched. The relevant slopes for the three isotopy moves are analyzed below.

\ref{I1} Since the equality $-[a_0^i,\dots,a_{j-1}^i]^{-1}-[a_0^l,\dots,a_k^l,n_{k+1}^l+1]^{-1}=1$ holds and $n_0^3$ is at least as much as the denominator of the two fractions, we can peel off from $U_l$ and $U_3$ as many negative slices that the slope of the torus $T$ we get around $F_i$ via \ref{L2} equals $[a_{k_i}^i,\dots,a_{j+1}^i]$ in the torus basis. Torus $T$ bounds a basic slice with the torus of slope $[a_{k_i}^i,\dots,a_{j+1}^i-1]$ which cuts off positive outermost slices from $U_i$. The torus of slope $[a_{k_1}^i,\dots,a_{j+2}^i]$ then gives continued fraction block with the torus $T$.

\ref{I2} These are essentially isotopies from \cite[Proposition 6.4]{GLS}. Peeling off a single (negative outermost) basic slice from $U_1$ and $U_2$, we obtain a torus of slope $0$, expressed in the background basis, around $F_3$. It corresponds to the slope $[a_{k_3}^3,\dots,a_1^3]$, and forms a glued-together basic slice with the outermost continued fraction block (with inner slope $[a_{k_3}^3,\dots,a_1^3-1]$), and hence, a continued fraction block with $[a_{k_3}^3,\dots,a_2^3]$. 

\ref{I3} In the proof of Proposition \ref{prop:ot}, for structures of the second kind (\ref{O2}), we have obtained that the slope of the torus $T$ built via \ref{L2} from the two tori of slopes $[a_0^1,\dots,a_J^1,3]^{-1}+1$ around $F_1$ and $[a_0^2,\dots,a_K^2,2]^{-1}$ around $F_2$, which sum up to zero, equals $-\frac{1}{D}$ for $D$ the denominator of $[a_0^1,\dots,a_J^1,3]^{-1}$. By assumption, $D$ equals $a_0^3$, moreover, $-\frac{1}{a_0^3}$ is in the torus basis expressed as $[a_{k_3}^3,\dots,a_2^3]$. Thus, torus $T$ bounds a basic slice with the torus of slope $[a_{k_3}^3,\dots,a_2^3-1]$ and a continued fraction block with the torus of slope $[a_{k_3}^3,\dots,a_3^3]$ in the slicing of $U_3$.
\endproof

\proof[Proof of Theorem \ref{thm}]
In Section \ref{Sec2} we have reduced Theorem \ref{thm} to Theorem \ref{thm'}.

Suppose we are given a contact surgery presentation $\xi$ as in Figure \ref{fig:SFS}, whose full path $\mathcal P_\xi$ properly ends. By Theorem \ref{thm:c} such a presentation describes a tight structure. Furthermore, all presentations which share this same path induce the same $\spinc$ structure (defining the same class in $\widehat{HF}(-M)$ \cite{OSz}) and have the same $3$-dimensional invariant (which is through the equality $d_3(\xi)=d(M,\mathbf t_\xi)$, Theorem \ref{thm:d=d3}, determined by $\mathbf t_\xi$). In Corollary \ref{cor:paths} we identify how presentations in the same path are related to each other, and in Proposition \ref{prop:isotopies} we realize all these relations by contact isotopies. Any further isotopies are, of course, excluded by the fact that different paths present non-homotopic bundles.

On the other hand, if the path $\mathcal P_\xi$ drops out (fails the tightness criterion), the corresponding structure $\xi$ admits one of the features recognized in Corollary \ref{cor:paths} or it can be walked through presentations, related by the above isotopy moves, to some presentation which admits such a feature. Finally, for these structures either the existence of a torus with critical slope or the Gluing Lemma argument proves their overtwistedness, in Proposition \ref{prop:ot}. 

This finishes the proof of Theorem \ref{thm'}, originally stated as Theorem \ref{thm}. 
\endproof


\end{document}